\documentclass[11pt, 3p]{elsarticle}

\usepackage{lineno}
\modulolinenumbers[5]

\usepackage[short]{optidef}
\usepackage[english]{babel}
\usepackage[utf8]{inputenc}
\usepackage{mathtools}
\usepackage{hyperref}
\usepackage{amsthm}
\usepackage{amssymb}
\usepackage{comment}
\usepackage{graphicx}
\usepackage{subcaption}
\usepackage{algpseudocode} 
\usepackage{algorithm}
\usepackage[table,xcdraw]{xcolor}
\usepackage{multirow}
\usepackage{rotating}
\usepackage{pdflscape}
\usepackage{longtable}
\usepackage{subcaption}
\usepackage{geometry}
\usepackage{url}
\usepackage{booktabs}
\usepackage{caption}
\usepackage{array}

\journal{}

\biboptions{authoryear}

\usepackage{setspace}
\begin{document}

\begin{frontmatter}

\title{Optimal Placement of Nature-Based Solutions for Urban Challenges}

\author{Diego Maria Pinto}
\ead{diegomaria.pinto@iasi.cnr.it}
\address{Institute for Systems Analysis and Computer Science ``Antonio Ruberti'', \\ National Research Council, Via dei Taurini 19, 00185 Rome, Italy}
\author{Davide Donato Russo\corref{cor1}}
\ead{davidedonato.russo@iasi.cnr.it}
\address{Institute for Systems Analysis and Computer Science ``Antonio Ruberti'', \\ National Research Council, Via dei Taurini 19, 00185 Rome, Italy}
\author{Antonio M. Sudoso}
\ead{antoniomaria.sudoso@uniroma1.it}
\address{Department of Computer, Control and Management Engineering ``Antonio Ruberti'', \\ Sapienza University of Rome, Via Ariosto 25, 00185 Rome, Italy \\ \& \\ Institute for Systems Analysis and Computer Science ``Antonio Ruberti'', \\ National Research Council, Via dei Taurini 19, 00185 Rome, Italy}

\cortext[cor1]{Corresponding author}

\begin{abstract}
Increased urbanization and climate change intensify urban heat islands and degrade air quality, making current mitigation strategies insufficient. Nature-based solutions (NBSs), such as parks, green walls, roofs, and street trees, offer a promising means to regulate urban temperatures and enhance air quality. However, determining their optimal placement to maximize environmental benefits remains a pressing challenge. Leveraging Operational Research (OR) tools, we propose a Mixed-Integer Linear Programming (MILP) model that integrates multiple factors, including urban challenges, physical constraints, clustering techniques, convolution theory, and fairness considerations. This model determines the optimal placement of NBSs by addressing metrics such as ground temperature, air quality, and accessibility to green spaces.
Through several case study analyses, we demonstrate the effectiveness of our approach in improving environmental and social indicators. This research holds implications for policy and practice, empowering urban planners and policymakers to make informed decisions regarding NBS implementation. Such decisions ensure that investments in urban greening yield maximum environmental, social, and economic benefits.
\end{abstract}


\end{frontmatter}

\newpage

\section{Introduction} \label{introduction}
Climate change and the current trends in urbanization make city resilience a clear priority \citep{valentina2019future, KIM2022665}. 
Nature-based Solutions (NBSs) are interventions through natural elements that simultaneously address social, economic, and environmental issues, thereby presenting a multifunctional, solution-oriented approach to increase sustainability \citep{dorst2019urban}.

\begin{figure}[!ht]
     \centering
     \begin{subfigure}[b]{0.3\textwidth}
         \centering
         \includegraphics[width=\textwidth]{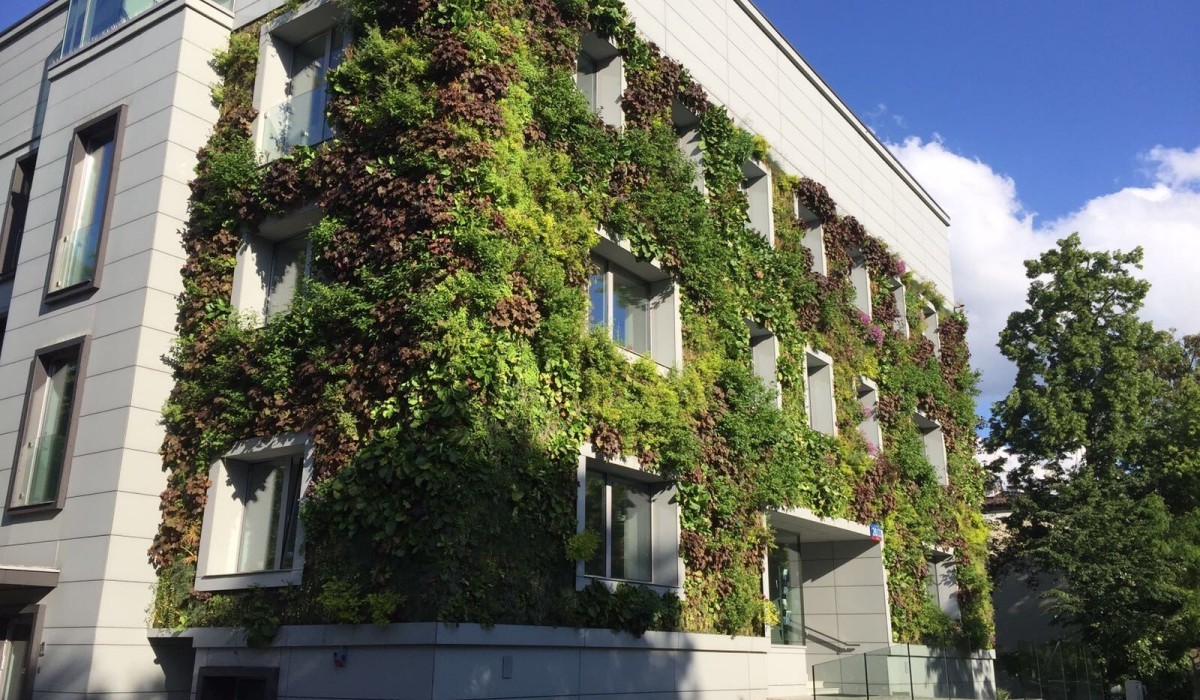}
        \caption{Green Wall in Warsaw (Poland)}
         \label{fig:Greenwall}
     \end{subfigure}
     \begin{subfigure}[b]{0.3\textwidth}
         \centering
         \includegraphics[width=\textwidth]{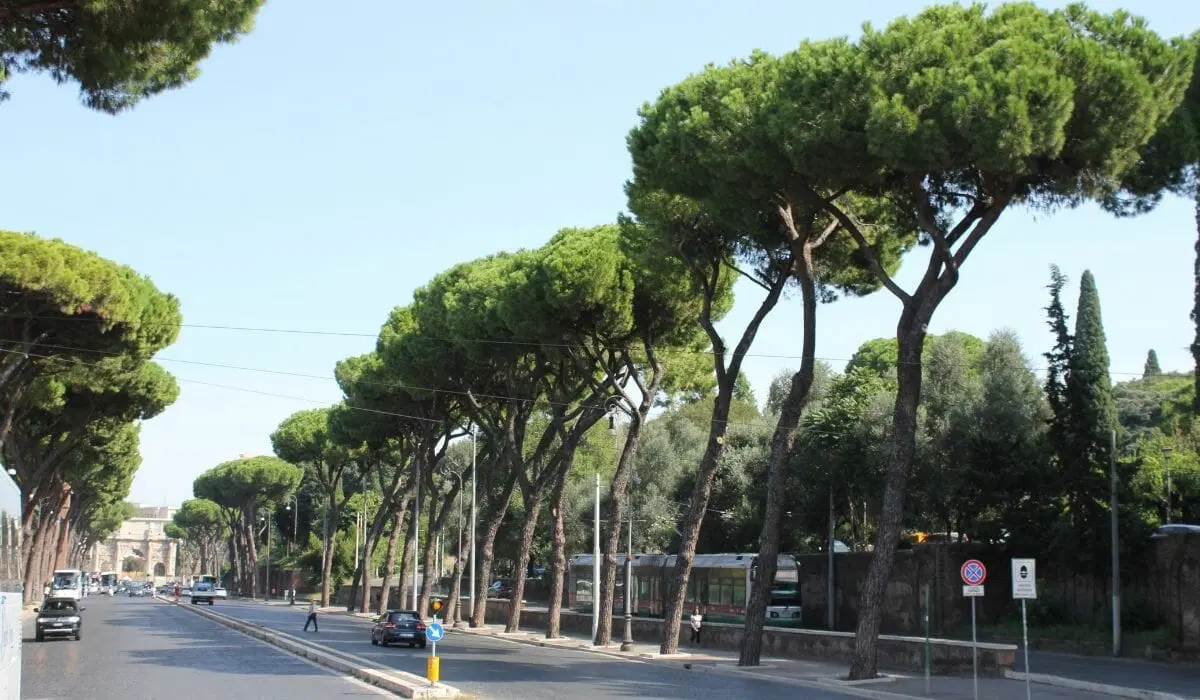}
        \caption{Street Trees in Rome (Italy)}
         \label{fig:Streettrees}
     \end{subfigure}
     \begin{subfigure}[b]{0.3\textwidth}
         \centering
         \includegraphics[width=\textwidth]{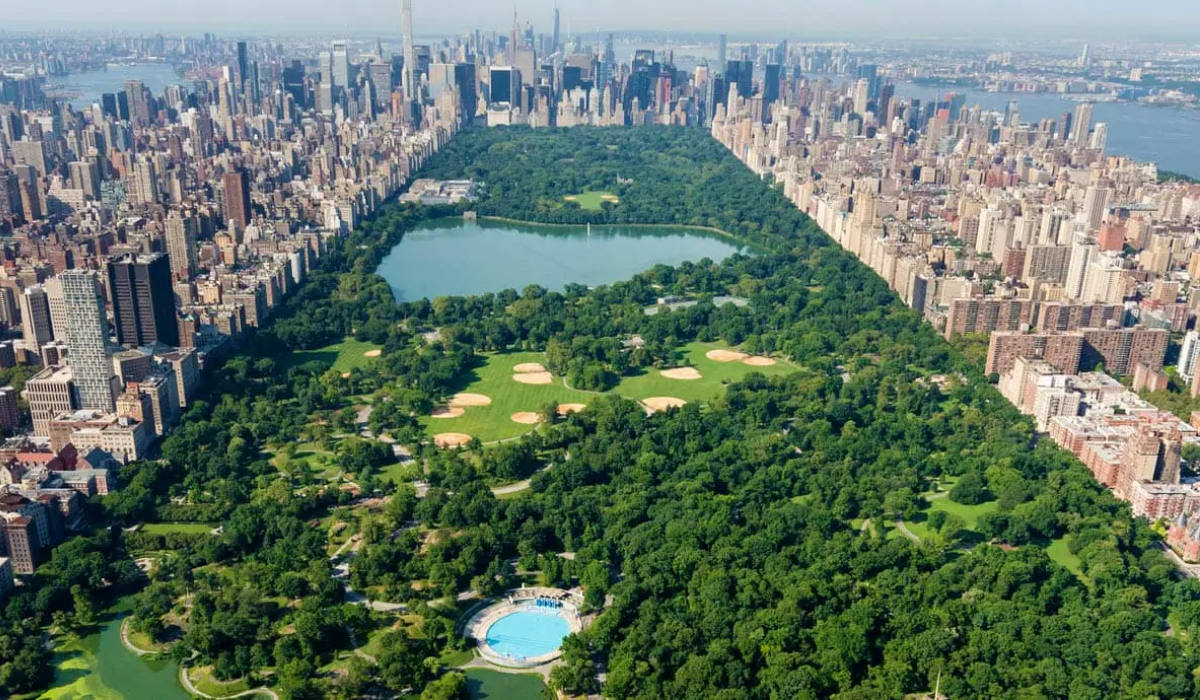}
        \caption{Urban Park in New York (USA)}
         \label{fig:Urbanpark}
     \end{subfigure}
        \caption{Three examples of NBSs.}
        \label{fig:NBSExamples}
\end{figure}

These solutions aim to create more resilient urban environments while enhancing the well-being of the residents \citep{HE2024122218}. Some key NBSs used in urban greening are urban parks, green spaces, green roofs and walls, tree planting, green corridors, permeable pavements, and urban farming. Such NBSs have been shown to offer numerous benefits, including improved air and water quality, reduced urban heat island effects, and enhanced biodiversity. For instance, installing vegetation on rooftops and vertical surfaces helps to reduce heat absorption, improve energy efficiency, and mitigate stormwater runoff; using permeable materials for sidewalks, roads, and parking lots allows rainwater to infiltrate the ground, reducing flooding and the burden on stormwater systems\citep{JUN20131003}; increasing the number of trees and establishing green corridors throughout the city enhances biodiversity, provides shade, and contributes to carbon sequestration. Examples of NBSs in urban greening are shown in Figure \ref{fig:NBSExamples}.

The literature addressing greening benefits in urban contexts is extensive, and quantification of greening effects usually requires relatively complex and integrated models \citep{leibowicz2020urban}. 
An European Union (EU)-wide assessment is presented by \cite{quaranta2021water}. This study quantifies the benefits of urban greening in terms of green water availability, reduction of cooling costs, and CO2 sequestration from the atmosphere across different climatic scenarios. Instead, the reader can refer to the work by \cite{rowe2011green} and references therein for a discussion on how green roofs influence air pollution, carbon dioxide emissions, and carbon sequestration. Additional references to greening benefits and their quantification in urban contexts can also be found in \citep{nowak2006air, currie2008estimates, gillner2015role, skoulika2014thermal, jamei2021review, oquendo2022systematic, segura2022street, susca2022effect}.


The implementation of NBSs can benefit from a data-driven and performance-based planning approach, which relies on empirical data and measurable outcomes to guide decision-making \citep{kabisch2017nature, dorst2019urban}. Unlike traditional planning methods, a performance-based approach focuses on achieving specific outcomes and can adapt to varying urban challenges, making it particularly well-suited for dynamic urban environments. Indeed, the importance of strategic planning in urban greening cannot be overstated \citep{rosenhead1981operational, ramshani2020optimal, DOUKAS20201}. Planning is essential not only for the integration of green spaces into the urban fabric but also for maximizing their expected ecological, social, and economic benefits. Therefore, effective urban greening should target the following key objectives: (i) \textit{Urban Heat Island Mitigation} - green spaces strategically placed in areas with high urban heat island effect can help cool down the city and mitigate the impacts of extreme heat on public health and energy consumption; (ii) \textit{Air Quality Improvement} - when locating green zones in areas with high air pollution or traffic congestion, the vegetation can act as a natural air filter, absorbing pollutants and improving air quality; (iii) \textit{Equitable Distribution} - strategic planning helps ensure that green spaces are fairly distributed across the city, providing access to nature and its benefits to all residents, regardless of their socioeconomic status or location \citep{morais2011evaluation}. Moreover, urban greening strategies are considered effective when they are supported by geo-referenced data for each specific Urban Challenge (UC) measure, such as ground temperature and air quality parameters. Equally important is the accurate quantification of the performance of NBSs, which involves assessing the tangible benefits these solutions deliver for each UC \citep{DIPIRRO20221, di2023cost}.

Considering the relevance of Operational Research (OR) in urban planning \citep{rosenhead1981operational, hua2019integer, DYE2024654}, in this paper we exploit OR tools to develop a flexible mathematical formulation to locate multiple NBSs while considering real-world constraints. Specifically, we make the following contributions: 

\begin{itemize} 
    \item We present a Mixed-Integer Linear Programming (MILP) model that optimizes the placement of various NBSs within urban environments, addressing metrics of different Urban Challenges (UCs), including ground temperature, air quality, and the population's accessibility to green areas.
    \item We evaluate the effectiveness of the MILP formulation using real-world case studies from Italian cities, selected for their diverse urban characteristics, climatic conditions, and geographic contexts.
    \item Data for both NBSs and UCs are gathered from state-of-the-art literature. This information is then integrated with national and international databases (e.g., population and land use) to create a set of benchmark instances. To foster further research, both the source code and datasets are made publicly available at \url{https://github.com/daviderusso/OptimalNBS}.
\end{itemize}

The remainder of the paper is organized as follows. Section \ref{litterature_review} provides a review of previous studies related to NBS localization. Section \ref{methodology} introduces the methodology underlying the formulation. Section \ref{formulation} presents the proposed MILP  model. Section \ref{experimental_setup} describes the experimental setup, including the instance creation process and parameters required to apply the model to real-world case scenarios. Section \ref{results} reports the experimental results and offers insights and an in-depth analysis of the resulting planning benefits. Finally, Section \ref{conclusions} presents conclusions and future research perspectives.

\section{Literature review} \label{litterature_review}
In this section, we review OR-based methods for urban greening planning. While there is a vast body of literature on the advantages of urban greening, few studies have applied OR tools to exploit those advantages within urban planning \citep{ehmke2016data}.

\cite{yoon2019multi} propose a multi-objective planning model to locate green spaces based on multiple effects, such as cooling and enhancement of connectivity. This optimization problem is solved using meta-heuristic algorithms and applied to synthetic landscapes. In \citep{yang2023multi}, a genetic algorithm is used to minimize flood risks according to the performance assessment of an urban drainage system model. In recent literature, urban (storm) water management is indeed frequently addressed by planning NBSs, as in \citep{saeedi2023planning} using genetic algorithms, in \citep{barah2021optimizing} with a two-stage stochastic programming model under precipitation uncertainty, and in \citep{Liu2023} where also rural contexts are considered. \cite{ramshani2020optimal} propose an approach to investigate the joint placement of photovoltaic (PV) panels and green roofs to improve the output efficiency of PV panels. They develop a two-stage stochastic programming model to incorporate PV panel/green roof placement decisions under different climate models to maximize the overall profit from energy generated and saved.

In \citep{castro2022optimizing}, social fairness is incorporated into the optimization process by combining social equity with hydro-environmental performance through a novel application of the Gini coefficient. Similarly, \citep{li2022spatial} addresses fairness in the planning of urban green spaces, aiming to minimize the conversion costs of newly added green areas while accounting for spatial equity constraints. \cite{xiao2002using} use an evolutionary algorithm to generate solutions to a site-search problem. This trend toward incorporating fairness in OR solutions continues indeed to gain importance \citep{boresta2024bridging}.

In \citep{zhang2017optimizing}, the authors present an optimization-based approach to reduce daytime and nighttime urban heat island effects. They propose a multi-objective model that maximizes the total sum of cooling benefits, either directly or indirectly. Thus, even if they consider cooling benefits only, the model is applied to evaluate the benefits achieved through simultaneous consideration of (i) direct cooling benefit for converting a certain area to green space, and (ii) indirect cooling benefit received from neighboring green spaces. Moreover, \cite{zhou2022optimization} address the heat island effect by optimizing the layout of urban green spaces using a genetic algorithm. Using a grid-based approach, \cite{hua2019integer} present an integer programming formulation as an urban design tool for residential projects. This tool considers both sunlight-gain rules through 3D modeling and routes connecting all building entrances to public transportation. \cite{pribadi2017optimizing} introduce mixed-integer nonlinear programs to define the spatial shape of urban green spaces. Their goal is to maximize the potential of these spaces to deliver various ecosystem services. For a recent literature review on the selection and placement of NBSs in the context of climate adaptation, readers can refer to the work by \cite{esmail2022greening} and \cite{capgras2023optimisation}.


\cite{viti2022knowledge} point out that current studies on NBSs often rely on isolated, case-specific approaches, with no shared, comprehensive methodology. This fragmented approach hinders the adoption of NBSs and fails to provide practical guidance for urban planners and city managers. They also suggest that future research should focus on applying systematic assessments of the benefits of NBSs. Building on this research gap, our work presents a novel contribution. Indeed, unlike previous studies that target single goals such as cooling or flood mitigation, our approach simultaneously considers a broader spectrum of NBSs and benefits over UCs. To the best of our knowledge, no prior research has integrated these diverse factors into a unified OR-based optimization framework.

\section{Methodology overview} \label{methodology}
\begin{figure}[!ht]
    \centering
    \includegraphics[scale=0.25]{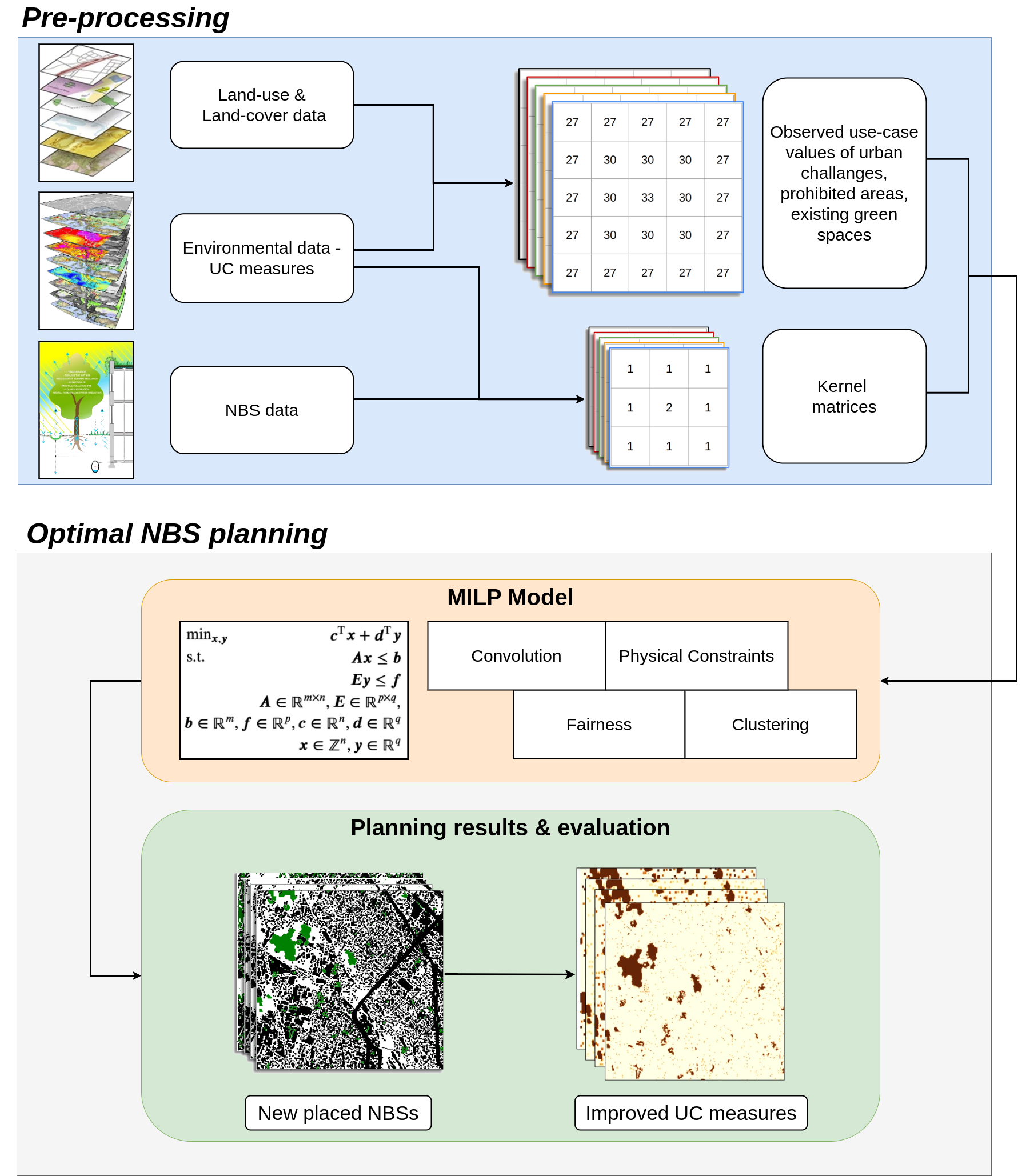}
    \caption{The fundamental components of our approach include pre-processing and optimal NBS planning phases. During the pre-processing phase, we collect all required data to stack a detailed set of matrices representing the present use case, alongside kernels for assessing the NBSs' impact. The subsequent phase involves the MILP formulation, which integrates convolution, fairness, clustering, and physical constraints. The output of this phase is a list of newly placed NBSs and the corresponding improvements in UC measures.}
    \label{fig:methodology_flowchart}
\end{figure}

The objective of our modeling approach is to locate various NBSs to enhance metrics of different Urban Challenges (UC), including ground temperature, air quality, and the population's accessibility to green areas. In the following, we set out the building blocks of the methodology, which is structured into two main phases: pre-processing and optimal NBS planning, as shown in Figure \ref{fig:methodology_flowchart}.

In the pre-processing phase, we collect data related to the area of interest. This includes (i) land-use and land-cover data; (ii) environmental data, such as ground temperature and air pollution; (iii) data on the effectiveness of NBSs in addressing UCs. Detailed information on the data collection and integration process can be found in Section \ref{experimental_setup}. We discretize the urban landscape into a fixed number of locations according to the spatial resolution of land use and land cover maps. This grid-based approach facilitates the localization of NBSs within the area of interest. Indeed, each cell within the obtained grid corresponds to a potential location for an NBS installation. 
During pre-processing, we also collect physical constraints, including lists of prohibited areas (e.g., rivers, streets, and other obstacles) and existing green spaces covering each cell. These data, combined with environmental data, are used to create layered maps that describe the area of interest.

\begin{figure}[!ht] \label{fig:kernel_example}
    \centering
        \begin{subfigure}[b]{0.23\textwidth}
        \centering
        \includegraphics[width=\textwidth]{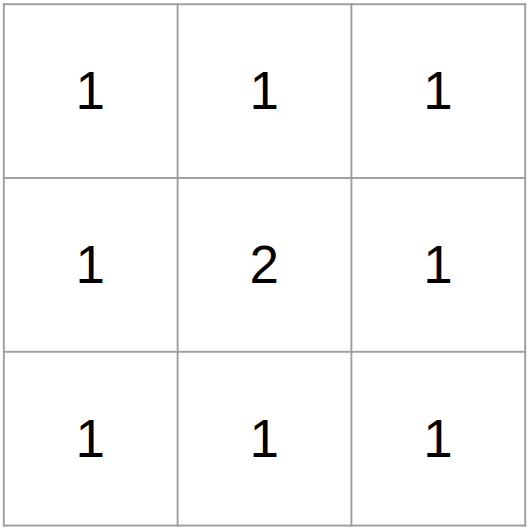}
        \caption{}
        \label{fig:kernel}
    \end{subfigure}
    \hfill
    \begin{subfigure}[b]{0.23\textwidth}
        \centering
        \includegraphics[width=\textwidth]{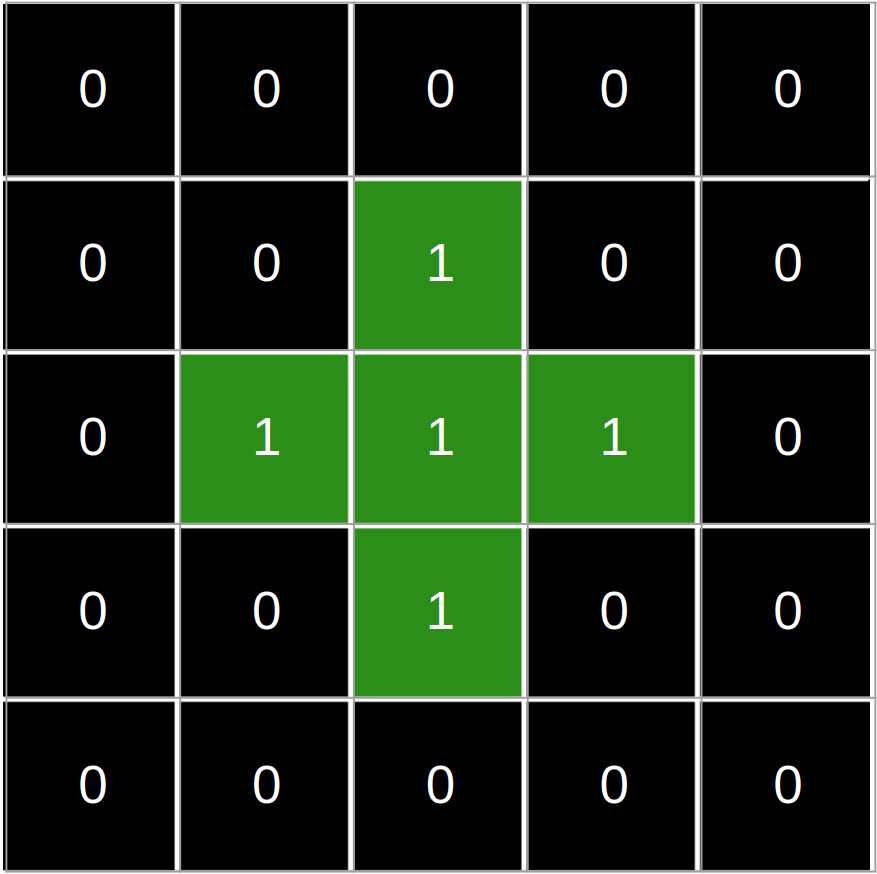}
        \caption{}
        \label{fig:newGreen_val}
    \end{subfigure}
    \hfill
    \begin{subfigure}[b]{0.23\textwidth}
        \centering
        \includegraphics[width=\textwidth]{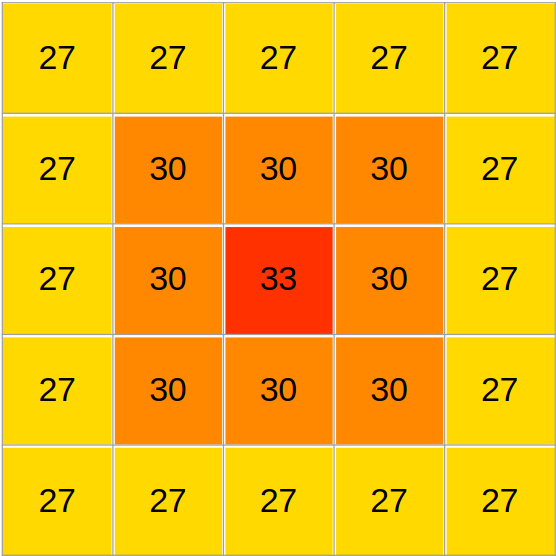}
        \caption{}
        \label{fig:original}
    \end{subfigure}
    \\
    \begin{subfigure}[b]{0.23\textwidth}
        \centering
        \includegraphics[width=\textwidth]{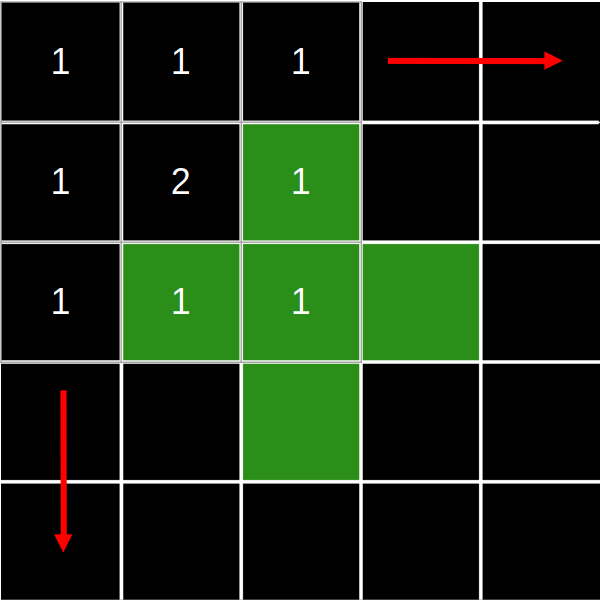}
        \caption{}
        \label{fig:Kernel_on_map}
    \end{subfigure}
    \hfill
    \begin{subfigure}[b]{0.23\textwidth}
        \centering
        \includegraphics[width=\textwidth]{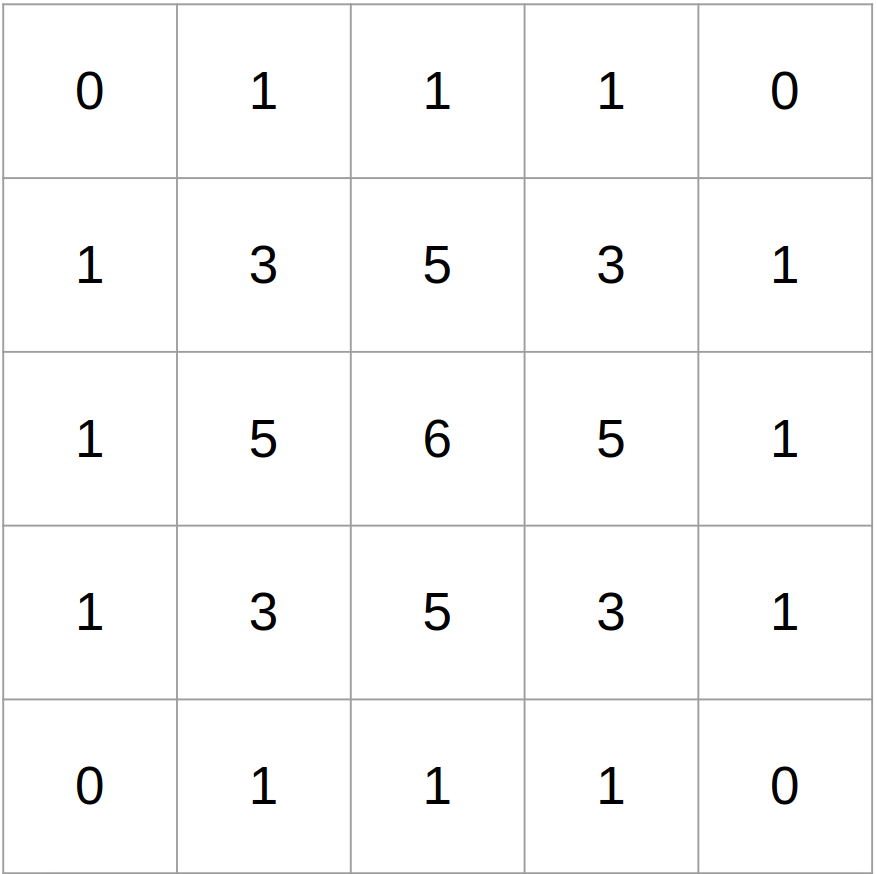}
        \caption{}
        \label{fig:reduction_value}
    \end{subfigure}
    \hfill
    \begin{subfigure}[b]{0.23\textwidth}
        \centering
        \includegraphics[width=\textwidth]{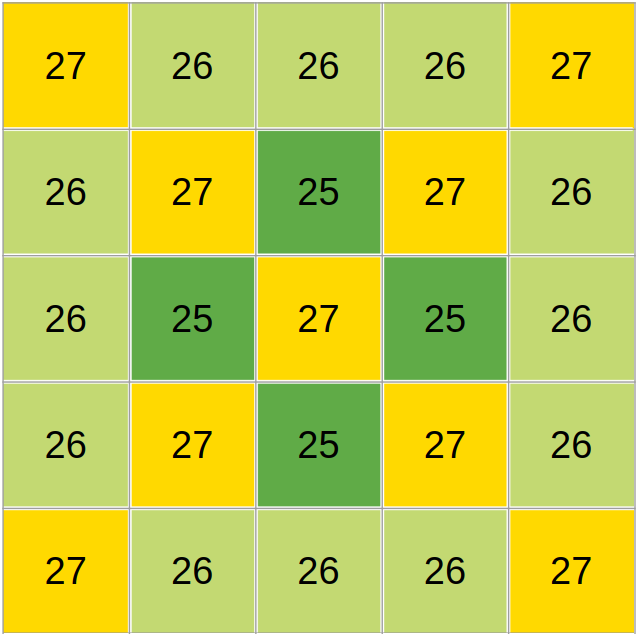}
        \caption{}
        \label{fig:results}
    \end{subfigure}
    \caption{The proposed convolution approach to assess the influence of NBSs on the urban grid. The kernel employed to assess the NBS influence on the map is depicted in Figure \ref{fig:kernel}. Figure \ref{fig:newGreen_val} represents a sample implementation of a new NBS (the green cells corresponding to the 1 entries). The kernel is applied to the newly installed NBSs as shown in Figure \ref{fig:Kernel_on_map}. This convolution is applied to evaluate the reduction of the current temperature whose heatmap and associated values are depicted in Figure \ref{fig:original}. As a result, the temperature decreases by the value depicted in Figure \ref{fig:reduction_value}. Figure \ref{fig:results} presents the updated heatmap obtained by subtracting the values in Figure \ref{fig:reduction_value} from that of the initial heatmap in Figure \ref{fig:original}. Before the NBS installation (Figure \ref{fig:original}), the highest peak temperature was $33^{\circ}$. After the installation (Figure \ref{fig:results}), it is reduced to $27^{\circ}$.} 
    \label{fig:KernelApplicationExample}
\end{figure}

Drawing inspiration from image filtering techniques, we estimate the effect of NBS installations through convolution operations \citep{capobianco2021image}. In image processing, convolution involves sliding a filter, referred to as a ``kernel matrix'', over an image. This process consists of performing element-wise multiplication between the kernel matrix and the corresponding region of the image, and then summing the results to determine the value of each output pixel in the filtered image \citep{XIAO2024124935}. Consequently, each output pixel represents a weighted sum of its neighboring input pixels. Similarly, we use a convolution approach to assess the influence of an NBS on the urban grid, where the kernel matrix defines the neighborhood impact of that NBS. By convolving this kernel with the grid of candidate locations, we can evaluate the overall effect of an NBS at each point on the grid. This procedure helps identify how each UC measure is affected by NBS interactions with adjacent cells, whether they host a different NBS type or the same one. To achieve this, we construct a kernel matrix for each NBS type and UC measure. More details regarding the kernel's construction can be found in Section \ref{experimental_setup}. As an example, Figure \ref{fig:KernelApplicationExample} illustrates the convolution between a $3 \times 3$ kernel (Figure \ref{fig:kernel}) and a $5 \times 5$ matrix (Figure \ref{fig:newGreen_val}). The kernel matrix models the reduction of the temperature values, as reported in the heatmap of Figure \ref{fig:original}. This reduction is computed considering the newly installed NBSs, represented by the green cells in Figure \ref{fig:newGreen_val} (i.e., the entries with value 1). The resulting temperature reduction is shown in Figure \ref{fig:results}. 

In the second phase, we employ a MILP model to determine the optimal locations for installing NBSs from a set of candidate sites. We do not directly address the sizing of NBSs, as each cell in the grid is constrained by its spatial resolution. However, to accommodate larger NBS types, such as urban parks, we apply clustering techniques to combine adjacent cells, allowing us to handle larger contiguous areas. This approach ensures our model can support different NBS sizes and configurations. The MILP model incorporates convolution operations to evaluate the spatial impact of NBSs on the selected UC measures. Additionally, we consider physical constraints, such as prohibited zones and existing green spaces, to ensure that the solution aligns with real-world scenarios. These constraints also account for the compatibility between NBS types and the land-use characteristics of each cell. Finally, guided by principles from the European Environment Agency\footnote{\url{https://climate-adapt.eea.europa.eu/}}, our model also considers fairness to ensure equitable access to green spaces across the area. Once the MILP model is solved, the results provide a list of recommended NBS locations along with the expected improvements in UC measures, helping planners make informed decisions.

\section{Mixed-integer linear programming formulation}\label{formulation}
In this section, we present the mathematical formulation of the MILP model.
We consider a rectangular grid where each cell has a spatial resolution of $s$ meters, i.e., every cell is a square of $s \times s$ m$^2$ on the ground. Let $G$ be the set of all possible locations for NBS installations within this grid. Specifically, $G$ consists of $W \times H$ pairs $(i, j)$ representing potential installation sites. Additionally, let $T$ be a set of different NBS that can be placed in each location $(i, j) \in G$. We define $c^t$ as the cost (e.g., installation and maintenance costs) for a single cell of NBS $t \in T$. Moreover, let $Budget$ represent the total available budget for NBS installation. A set of forbidden locations $(i, j)$ can be defined within $G$, namely $F^t \subseteq G$, containing $(i, j) \in G$ where NBS type $t \in T$ cannot be installed due to the presence of physical constraints (i.e., rivers, buildings, streets). Similarly, $E^t \subseteq G$ denotes the set of locations $(i, j) \in G$ where an NBS type $t \in T$ has already been installed.
As discussed in Sections \ref{introduction} and \ref{litterature_review}, different NBSs have impacts on UC measures. Accordingly, we define the set $U$ of UC measures (e.g., ground temperature and air pollution) that can be tackled through NBSs by exploiting their beneficial effects, like those linked to air quality and urban heat islands (e.g., reduction of particulate matter (PM), mitigation of temperature, reduction of noise pollution). We denote by $A^{u}$ the $W \times H$ matrix of observed values $a^{u}_{ij}$ of $u \in U$ (e.g., observed ground temperature, observed PM concentrations). For each $u \in U$, these values may be related to specific sampling times or long-term averages. As stated in the literature, the maximum benefit achievable by installing an NBS in a specific area is constrained by physical and practical limitations \citep{seddon2020understanding}. For this reason, we denote by $\delta^u$ the maximum achievable reduction of $u \in U$ for the whole area of interest, given the current observed values $A^{u}$.

As discussed in Section \ref{methodology}, we use convolution operations to evaluate the effect of NBS installations. Indeed, for each type of NBS $t \in T$ and each UC measure $u \in U$, we define the corresponding kernel matrix $K^{ut}$ of size $\bar{w}^{ut} \times \bar{h}^{ut}$. Note that, the kernel sizes may vary based on the characteristics of each $t \in T$ and $u \in U$. The specific details on how to model convolution operations are presented in the mathematical formulation within this section. The adopted notation, including parameters and indices, is outlined below and summarized for quick reference in Table \ref{tab:Notation}.

\begin{table}[!ht]
    \centering
    \footnotesize
    \begin{tabular}{l|l}
    \toprule
    \textbf{Sets and indices} & \textbf{Description} \\
    \midrule
     ${T}$    & Set of Nature-Based Solutions (NBSs) \\
     ${U}$    & Set of Urban Challenge (UC) measures \\ 
     $G$      & Set of candidate locations for NBSs \\
     $F^t$    & Set of forbidden locations for the NBS of type $t \in {T}$\\
     $E^t$    & Set of locations containing pre-existing NBS of type $t \in {T}$\\
     $t$      & Index for the NBS of type $t \in T$ \\
     $u$      & Index for the UC measure $u \in U$ \\
     $(i, j)$ & Pair of indices for the location $(i, j) \in G$ \\
     $Q^t$    & Set of cluster indices for NBS of type $t \in T$\\
     $P^t_q$  & Set of $(i,j) \in G$ assigned to the cluster $q \in Q^t$\\   
     \midrule
     \textbf{Parameters} &  \textbf{Description} \\
     \midrule
     $s$             & Spatial resolution of a grid cell, expressed in meters\\
     $W$             & Width of the area of interest divided by $s$ meters \\
     $H$             & Height of the area of interest divided by $s$ meters \\
     $c^t$           & Installation cost for the NBS $t \in T$ \\
     $A^u$           & $W \times H$ matrix of observed values of the UC measure $u \in U$ \\
     $a^{u}_{ij}$    & Current observed value for $u \in U$ in $(i, j) \in G$\\
     $\delta^u$      & Maximum reduction value of $u \in U$ \\    
     $K^{ut}$        & Kernel matrix of NBS $t \in T$ for UC measure $u \in U$ \\
     $\theta^u_\textrm{max}$ &  Relative importance of the highest values of $u \in U$ in the objective function \\  
     $\theta^u_\textrm{avg}$ & Relative importance of the average values of $u \in U$ in the objective function \\ 
     $\theta^c$      & Relative importance of the NBSs costs in the objective function \\ 
     $\theta^f$      & Relative importance of the fairness term in the objective function \\
     $Budget$        & Total available budget \\
     $M^u$           & Big-M value for $u \in U$\\
     $\bar{w}^{ut}$  & Kernel width for $u \in U$ and $t \in T$ normalized by $s$ meters\\
     $\bar{h}^{ut}$  & Kernel height for $u \in U$ and $t \in T$ normalized by $s$ meters\\
     $\Pi_{ij} $   & Population percentage in $(i, j) \in G$\\
     \midrule
     \textbf{Variables} &  \textbf{Description} \\
     \midrule 
     $x_{ij}^{t} \in \{0, 1\}$               & Indicating if an NBS type $t \in T$ is used in the tile $(i, j) \in G$ or not \\
     $y^u_{ij} \in \{0, 1\}$                 & Indicating if the value of $z^u_{ij}$ exceeds the value of $\delta^u$ \\
     $\lambda^t_q \in \{0, 1\}$                      & Indicating whether the cluster $q \in Q^t$ is used or not for an NBS $t \in T$ \\
    $z_{ij}^{u} \in \mathbb{R}^+$           & The result of the application of the kernel in $(i, j) \in G$ for $u \in U$\\
     $\bar{z}^{u}_{ij} \in \mathbb{R}^+$     & The actual value of the reduction provided by all NBS to $u \in U$ in tile $(i,j) \in G$ \\
     $z_{\textrm{max}}^{u} \in \mathbb{R}^+$ & The peak values for the UC measure $u \in U$\\
     $z_{\textrm{avg}}^{u} \in \mathbb{R}^+$ & The average value across all locations $(i,j) \in G$ of UC measure $u \in U$\\
     {$f_{ij} \in \mathbb{R}^+$}               & {Product between $\Pi_{ij}$ and the application of the fairness kernel in $(i, j) \in G$ for all $t \in T$}\\
     \bottomrule
    \end{tabular}
    \caption{Summary of the notation adopted for the formulation, including parameters, indices, sets, and variables.}
    \label{tab:Notation}
\end{table}

The model considers the following variables: $x_{ij}^t \in \{0, 1\}$ where $x_{ij}^t = 1$ if an NBS type $t \in T$ is used in the location $(i, j) \in G$ and $x_{ij}^t = 0$ otherwise; $z_{ij}^{u} \in \mathbb{R}^+$ is the result of the application of the kernel for all the NBS type $t \in T$ over the UC measure $u \in U$ for the location $(i, j) \in G$; $\bar{z}^{u}_{ij} = \min(z_{ij}^{u}, \delta^{u})$ is the actual value of the reduction provided by all NBS to $u \in U$ in location $(i,j) \in G$; $z_{\textrm{max}}^{u} \in \mathbb{R}^+$ represents the peak values for the UC measure $u \in U$; $z_{\textrm{avg}}^{u} \in \mathbb{R}^+$ represents the average across all the locations $(i,j) \in G$ of the UC measure $u \in U$. The MILP model is formulated as follows:
\begin{mini!}[2]
{}{\sum_{u \in U} \theta_{\textrm{max}}^{u} z_{\textrm{max}}^u + \sum_{u \in U}
   \theta_{\textrm{avg}}^{u} z_{\textrm{avg}}^u  + \theta^c \sum_{t \in T} \sum_{(i, j) \in G} c^t x_{ij}^t - \theta^f \sum_{(i, j) \in G} f_{ij}}
    {\label{prob:optimization_problem}}{}
\addConstraint{\sum_{t \in T} x_{ij}^t}{\leq 1}{\quad \forall (i, j) \in G \label{only_one_type}}
\addConstraint{\sum_{t \in T} \sum_{(i,j) \in G \setminus E^t} x_{ij}^t c^t}{\leq Budget \label{max_budget}}{}
\addConstraint{x_{ij}^t}{= 0}{\quad \forall (i, j) \in F^t : t \in T \label{forbidden_tile}}
\addConstraint{ x_{ij}^t}{= 1}{\quad \forall (i, j) \in E^t : t \in T \label{pre_existing_tile}}
\addConstraint{x_{ij}^t}{ = \lambda^t_q}{\quad \forall  t \in T, \ \forall q \in Q^t, \ \forall (i, j) \in P_q^t \label{constr:clustering}}
\addConstraint{z_{ij}^{u}}{= \sum_{t \in T} 
\sum\limits_{{i'}=1}^{\bar{w}^{ut}} \sum\limits_{{j'}=1}^{\bar{h}^{ut}} x^t_{I J} \times K^{ut}_{{i'}{j'}} }{\quad \forall (i, j) \in G, \ \forall u \in U : (I, J) \notin E^t\label{pixelComputEquation}}
\addConstraint{z_{ij}^u}{\le \delta^u + M^u (1 - y_{ij}^u)}{\quad \forall (i, j) \in G, \ \forall u \in U \label{linear1}}
\addConstraint{z_{ij}^u}{\ge \delta^u - M^u y_{ij}^u}{\quad \forall (i, j) \in G, \ \forall u \in U \label{linear2}}
\addConstraint{\bar{z}_{ij}^u}{\le z_{ij}^u}{\quad \forall (i, j) \in G, \ \forall u \in U \label{linear3}}
\addConstraint{\bar{z}_{ij}^u}{\le \delta^u}{\quad \forall (i, j) \in G, \ \forall u \in U \label{linear4}}
\addConstraint{\bar{z}_{ij}^u}{\ge z_{ij}^u - M^u (1 - y_{ij}^u)}{\quad \forall (i, j) \in G, \ \forall u \in U \label{linear5}}
\addConstraint{\bar{z}_{ij}^u}{\ge \delta^u - M^u y_{ij}^u}{\quad \forall (i, j) \in G, \ \forall u \in U \label{linear6}}
\addConstraint{z_{\textrm{max}}^u}{\ge a_{ij}^u - \bar{z}_{ij}^{u} }{\quad \forall (i,j) \in G, \ \forall u \in U \label{z_max_val}}
\addConstraint{z_{\textrm{avg}}^u}{= \frac{1}{W \times H}\sum_{(i, j) \in G} \left(a_{ij}^u - \bar{z}_{ij}^{u}\right)}{\quad \forall u \in U \label{z_avg_val}}
\addConstraint{f_{ij}}{= \Pi_{ij} \sum_{t \in T} \sum\limits_{{i'}=1}^{\bar{w}^{ft}} \sum\limits_{{j'}=1}^{\bar{h}^{ft}} x^t_{IJ} \times K^{ft}_{{i'}{j'}}}{\quad \forall (i, j) \in G  \label{fairness_calc}}
\addConstraint{x_{ij}^t}{\in \{0, 1\}}{\quad \forall (i, j) \in G, \ \forall t \in T}
\addConstraint{z_{ij}^u, \ \bar{z}_{ij}^u}{\in \mathbb{R}^+}{\quad \forall (i, j) \in G, \ \forall u \in U}
\addConstraint{z_{\textrm{max}}^u, \ z_{\textrm{avg}}^u}{\in \mathbb{R}^+}{\quad \forall u \in U}
\addConstraint{y_{ij}^u}{\in \{0, 1\}}{\quad \forall (i, j) \in G, \ u \in U}
\addConstraint{\lambda^t_q }{\in \{0, 1\}}{\quad \forall t \in T, \ q \in Q^t}
\addConstraint{f_{ij}}{\in \mathbb{R}^+}{\quad \forall (i, j) \in G}.
\end{mini!}

The objective function is given by the sum of four contributions. The first two terms aim to minimize UC measures within the area of interest, such as temperature and air pollution. Specifically, the first term focuses on reducing the peak values associated with each UC measure $u \in U$. We incorporate the average reduction value into the objective function to ensure that the model addresses urban challenges throughout the map, even when the peaks are no longer further improvable. This value is influenced by the installation of NBS across the entire area of interest, rather than solely in the proximity of the peaks associated with each urban challenge. Thus, the second term minimizes the average UC measure $u \in U$ across all locations $(i, j) \in G$. The third term minimizes the total cost of installing new NBSs, including both installation and long-term maintenance expenses. Finally, the last term promotes equitable accessibility to NBSs by maximizing the fairness in their spatial distribution across the area of interest. We point out that the order of magnitude of the terms in the objective function may differ considerably. For this reason, all the quantities in the objective function have been normalized in the range $[0, 1]$ to make them comparable. Moreover, for all $u \in U$, we introduce weights $\theta_{\textrm{max}}^{u}$, $\theta_{\textrm{avg}}^{u}$, representing the relative importance of each UC in the objective function, distinguished between the peaks and the average values of each $u \in U$. Furthermore, we use $\theta^c$ and $\theta^f$ as weights for the NBS costs and fairness terms. All these weights are non-negative and sum up to $1$.

The formulation includes the following constraints. Constraints \eqref{only_one_type} ensure that in each location $(i, j) \in G$ at most an NBS type $t \in T$ can be used; constraint \eqref{max_budget} imposes that the total cost must be within the available budget. Note that the locations in $E^t$ do not count for the total cost. Constraints \eqref{forbidden_tile} set to $0$ the variables that refer to locations belonging to the forbidden set for each NBS $t \in T$. As well, constraints \eqref{pre_existing_tile} set to $1$ the variables referring to locations of pre-existing NBSs.
Considering a UC measure $u \in U$, the effect over the location $(i,j) \in G$ of the whole set of NBSs is computed by constraints \eqref{pixelComputEquation}. Specifically, for each $(i,j) \in G$ and $u \in U$, the sum over $t \in T$ of the convolution between $K^{ut}$ and the matrix obtained by collecting variables $x^t_{IJ}$, where $I = i-\lfloor \frac{\bar{w}^{ut}}{2} \rfloor + {i'} -1$ and $J = j-\lfloor \frac{\bar{h}^{ut}}{2} \rfloor + j' - 1$ is computed. Then, constraints \eqref{pixelComputEquation} allow the computation of the kernel impact of newly installed NBS only, meaning that the convolution is computed if and only if $(I, J) \notin E^{t}$.

For each $u \in U$ and each location $(i, j) \in G$, the reduction of value is bounded between $0$ and $\delta^u$, that is
\begin{equation}
\label{eq:if_reduction}
    \bar{z}^u_{ij} =
    \begin{cases}
      z_{ij}^u, & \text{if}\ z_{ij}^u \leq \delta^u \\
      \delta^u, & \text{otherwise}.
    \end{cases}
\end{equation}
To model \eqref{eq:if_reduction}, we introduce a binary variable $y_{ij}^u$  such that $y_{ij}^u = 1$ if $z_{ij}^u \leq \delta^u$, and $y_{ij}^u = 0$ otherwise. Then, using $y_{ij}^u$, the value if $\bar{z}^u_{ij}$ can be expressed with the following nonlinear constraints: 
\begin{align}  
\label{constr:nonlinear_reduction}
    \bar{z}^u_{ij} = y_{ij}^u (z_{ij}^u - \delta^u) + \delta^u, \qquad \forall (i, j) \in G, \ \forall u \in U.
\end{align}
To linearize \eqref{constr:nonlinear_reduction}, we use a standard technique from the literature \citep{liberti2004reformulation}. Specifically, we introduce constraints \eqref{linear1} - \eqref{linear2} - \eqref{linear3} - \eqref{linear4} - \eqref{linear5} - \eqref{linear6}, where $M^u$ is a sufficiently large constant (big-M).
Constraints  \eqref{linear1} and \eqref{linear2} ensure that the binary variable $y_{ij}^u$ behaves correctly. Constraints \eqref{linear3} and \eqref{linear4} make sure that $\bar{z}_{ij}^u$ cannot exceed $z_{ij}^u$ and $\delta^u$, respectively. Finally, constraints \eqref{linear5} and \eqref{linear6} ensure that $\bar{z}_{ij}^u$ takes the appropriate value depending on the value of $y_{ij}^u$.
To estimate valid big-M quantities $M^u$ for each $u \in U$, these are set to the highest possible value achieved by installing the NBS $t \in T$ that offers the greatest improvement for $u$ across the area covered by the associated kernel.

Recall that the objective function's first term focuses on reducing the peak values associated with each UC measure $u \in U$. These peaks are represented by the variables $z_{\textrm{max}}^u$. In constraints \eqref{z_max_val}, these variables are set to be greater than or equal to the difference between $a_{ij}^{u}$ and $\bar{z}_{ij}^{u}$ over all the locations $(i, j) \in G$. Furthermore, the objective function's second term minimizes the average UC measure $u \in U$ across all tiles $(i, j) \in G$. These values are encoded using variables $z_{\textrm{avg}}^u$.  In constraints \eqref{z_avg_val}, $z_{\textrm{avg}}^u$ are set to be equal to the average difference between $a_{ij}^{u}$ and $\bar{z}_{ij}^{u}$, across all $(i, j) \in G$ of the UC measure $u \in U$.

As discussed in Section \ref{methodology}, we consider clustering constraints to represent contiguous areas suitable for NBS installation. Formally, for each $t \in T$ we denote by $Q^t$ the set of cluster indices. Furthermore, for each NBS $t \in T$, let $P^t_q = \{(i, j) \in G : (i,j) \ \textrm{is assigned to cluster} \ q \in Q^t\}$. Note that $\cup_{q \in Q^t}^{t \in T} P_q^t = G$ and $P_{q_1}^t \cap P_{q_2}^t = \emptyset$ for all $q_1 \neq q_2$ and $t \in T$. This assignment is determined in the pre-processing phase, as part of the instance creation process, detailed in Section \ref{experimental_setup}. In the model, variable $\lambda_q^t \in \{0, 1\}$ indicates whether the cluster $q$ is used ($\lambda^t_q = 1$) or not ($\lambda^t_q = 0$) for $t \in T$. Thus, constraints \eqref{constr:clustering} guarantee that all the locations that belong to the same cluster are simultaneously used for installing an NBS. Specifically, these constraints ensure that if $x_{ij}^t = 1$ for a given $(i, j) \in P^t_q$ then $\lambda^t_q = 1$. Therefore, variables $x_{ij}^t$ for all $(i, j) \in P^t_q$ are forced to be 1. If none of the variables $x_{ij}^t$ for a given $(i, j) \in P^t_q$ is equal to 0 then $\lambda^t_q = 0$. 


To keep into account fairness in the localization process, we model the accessibility of an NBS as an urban challenge. Specifically, for each NBS, we construct a kernel to assess the effect that an NBS installation has on the resident population. For each location \((i, j) \in G\), we define \(f_{ij}\) as the product of the resident population percentage \(\Pi_{ij}\) and the convolution of the kernel \(K^{tf}\) for each \(t \in T\) with the matrix formed by the variables \(x_{IJ}^t\). The indices \(I\) and \(J\) are calculated as $I = i-\lfloor \frac{\bar{w}^{ft}}{2} \rfloor + {i'} -1$ and $J = j-\lfloor \frac{\bar{h}^{ft}}{2} \rfloor + j' -1$. Note that $\bar{w}^{ft}$ and $\bar{h}^{ft}$ correspond to the dimensions of the fairness kernel, hence the superscript $f$. Constraints \eqref{fairness_calc} compute the fairness value for each location $(i, j) \in G$. The fairness term in the objective function is normalized between 0 and 1. Specifically, we perform the min-max normalization where the minimum value is given by the initial solution with pre-existing NBS, and the maximum is given by the hypothetical solution with an unlimited budget, thus installing only the NBS that provides the maximum improvement to fairness.

Problem \eqref{prob:optimization_problem} can be solved to optimality by using off-the-shelf MILP solvers. The computational effort grows with the number of NBSs, UC measures, and most importantly, with the grid dimension $W \times H$ which, in turn, depends on the spatial resolution $s$. As a consequence, using higher-resolution maps provides finer and more precise locations but also increases the computational burden associated with solving the optimization problem.

\section{Experimental Setup}\label{experimental_setup}

In this section, we describe the adopted experimental setup. Next, we provide details on how we gathered data used for the instance creation process. 

We selected four different NBSs, all dealing with green infrastructures. These are Green Roof (GR), Green Wall (GW), Street Tree (ST), and Urban Park (UP). As UC measures, we considered $\textrm{PM}_{10}$ concentration, $\textrm{PM}_{2.5}$ concentration, maximum ground temperature (TempMax), minimum ground temperature (TempMin), and the fair access to urban green areas (Fairness). All these UCs and NBSs are intriguing and interesting for this research, nonetheless, as one of the main strengths of this methodology, they can be adapted and personalized according to the characteristics of the test case.
In Table \ref{tab:NBS-Cost} we report the average installation (Inst. cost) and maintenance costs (Maint. cost) of the different NBS \citep{di2022facing, di2023cost}. Note that, the total cost (Tot. cost) is obtained by considering the sum of the maintenance cost plus the investment cost with a depreciation in 7 years. 

\begin{table}[H]
    \centering
    \footnotesize
    \begin{tabular}{|c|c|c|c|c|}
    \hline
        \textbf{NBS type} &  \textbf{Inst. cost (EUR/$\textrm{m}^2$)} &  \textbf{Maint. cost (EUR/$\textrm{m}^2$/yr)} & \textbf{Tot. cost (EUR/$\textrm{m}^2$/yr)} \\
        \hline
         Green Wall & 470 & 11.8 & 78.9 \\
         Green Roof & 310 & 7.8 & 52.0 \\
         Street Trees & 125 & 3.1 & 21.0 \\
         Urban Park & 225 & 5.6 & 37.8 \\ \hline
    \end{tabular}
    \caption{NBS cost estimation proposed by \cite{di2022facing, di2023cost}. The costs are reported in Euro for each $m^2$.}
    \label{tab:NBS-Cost}
\end{table}

The instance creation procedure involves several steps, from gathering information on UCs and NBSs, to defining the area for each test case. From the Italian territory, we identified the following cities: Bologna, Catania, Milan, Naples, Rome, and Turin. All of them have specific characteristics that make them suitable to be considered as test cases for our approach. These cities experience considerable environmental pollution, affecting both air and ground, which heightens their suitability for implementing NBSs. We gathered the following data from the Italian territory:
\begin{itemize}
    \item \textit{Land-Use and Land-Cover Data.} We retrieved land-use and land-cover data from \cite{landuse} to delineate areas suitable for locating NBS and areas where such placement is prohibited. This dataset, updated to 2023, provides a spatial resolution of 10 m.

    \item \textit{Air Quality Data.} Concentration data for $\textrm{PM}_{2.5}$, $\textrm{PM}_{10}$ from the European Environment Agency (EEA) dataset, updated to 2019. Air pollutant concentrations, with a spatial resolution of 1 km, are expressed in units of $\frac{\mu g}{m^3}$.

    \item \textit{Temperature Data.} We acquired maximum and minimum temperature data from \cite{DIPIRRO20221}. This dataset, updated to 2021, with a spatial resolution of $1$ km, includes daily temperature records. We used the data from July 19th as a benchmark. On this day, the maximum temperature span goes from a minimum of $4.95$ °C and a maximum of $35.60$ °C; the minimum temperature varies between -$2.29$ °C and $25.19$ °C.

    \item \textit{National Boundary Data.} To better characterize the area of interest, we utilized boundary information from the `Istituto Nazionale di Statistica' (ISTAT), updated to 2021.
    
    \item \textit{Population Density Data.} To provide a detailed evaluation of NBS distribution fairness relative to population density, we retrieved population density per $\textrm{km}^2$ data from ISTAT, updated to 2021.
\end{itemize}

Given the characteristics of the available datasets, the spatial resolution $s$ for the instance creation process and the MILP model is set to $10$ meters. Further details about the instances and their creation procedure are reported in Section \ref{results} and \ref{appendix:instance_details}, respectively.

For each NBS, its effect on the selected UCs must be properly identified. Therefore, information on their expected and proven impact must be retrieved from the relevant literature and data sources. We emphasize that conducting an exhaustive analysis of the effects of NBSs on various UCs is beyond the scope of this current research. We also acknowledge that a thorough understanding of these effects requires specific data collection and analysis tailored to each unique use case, whether it relates to a particular city, ecosystem, or type of NBS, and conducted over a defined period. 

Still, to provide an estimation of an NBS impact and populate the kernel entries, we averaged the target values presented in the reviewed literature \citep{nowak2018air, anderson2022lowering}. Overall, we created a total of 20 kernel matrices, one for each combination of four NBS and five UC. Furthermore, $\delta^u$ (the maximum reduction value for $u \in U$ within the area of interest) is computed as $\delta^u = 0.2 \times \max_{(i, j) \in G}(a_{ij}^u)$ for each $u \in U$. Note that, this estimation can be refined further, for example, by incorporating expert knowledge. A detailed explanation of the kernel entries can be found in \ref{appendix:kernel_details}, while the resulting kernel sizes and values are summarized in Table \ref{tab:kernel-info}. 


\begin{table}[!ht]
\centering
\footnotesize
\begin{tabular}{|c|c|c|c|c|c|c|}
\cline{3-7}
\multicolumn{2}{c|}{}& \multicolumn{5}{c|}{\textbf{UC measure}} \\ 
\hline
\multicolumn{2}{|c|}{\textbf{NBS type}} & \textbf{TempMax} & \textbf{TempMin} & \textbf{$\textrm{PM}_{2.5}$} & \textbf{$\textrm{PM}_{10}$} & \textbf{Fairness}  \\ \hline
\multirow{2}{*}{{Green Wall}}  & {Size}  & 5 $\times$ 5                & 3 $\times$ 3                & 5 $\times$ 5           & 5 $\times$ 5          & 5 $\times$ 5                 \\
\multicolumn{1}{|c|}{\textbf{}}           & {Range} & [0.10, 2.70]& [0.10, 1.90]& [0.10, 5.03]& [0.10 , 12.90]& [2.0, 6.0]\\ \hline
\multirow{2}{*}{{Green Roof}}  & {Size}  & 5 $\times$ 5                & 3 $\times$ 3                & 5 $\times$ 5           & 5 $\times$ 5          & 1 $\times$ 1                 \\ 
\multicolumn{1}{|c|}{\textbf{}}           & {Range} & [0.10, 2.00]& [0.10, 1.40]& [0.10, 2.51]& [0.10, 6.45]& [0.1, 2.0]\\ \hline
\multirow{2}{*}{{Street Tree}} & {Size}  & 5 $\times$ 5                & 3 $\times$ 3                & 3 $\times$ 3           & 3 $\times$ 3          & 3 $\times$ 3                 \\
\multicolumn{1}{|c|}{\textbf{}}           & {Range} & [0.10, 1.30]& [0.10, 0.70]& [0.10,  4.02]& [0.10, 10.32]& [0.1, 4.0]\\ \hline
\multirow{2}{*}{{Urban Park}}  & {Size}  & 5 $\times$ 5                & 3 $\times$ 3                & 7 $\times$ 7           & 7 $\times$ 7        & 11 $\times$ 11               \\ 
\multicolumn{1}{|c|}{\textbf{}}           & {Range} & [0.10, 3.50]& [0.10,  2.50]& [0.10,  5.03]& [0.10, 12.90]& [4.0, 10.0]\\ 
\hline
\end{tabular}%
\caption{Sizes and ranges for each kernel. Ranges are considered with the same measurement unit of the related UC measure.}
\label{tab:kernel-info}
\end{table}

We used the connected-component labeling method \citep{he2009fast} to automatically detect clusters of contiguous cells that are appropriate for installing structured NBSs. 
This method, also known as region labeling or region extraction, is widely used in computer vision to identify and label connected regions (components) in a binary matrix.
Specifically, we extract a binary matrix $B^t$ of size $W \times H$ from the set $G$ such that $B^t_{ij} = 1$ if the location $(i, j) \in G$ is suitable for installing an NBS of type $t \in T$, and 0 otherwise. Next, we identify connected regions of 1s and assign each a unique cluster index. Moreover, we extract regions with sizes ranging from a minimum of $5$ to a maximum of $50$ cells. 
In our experimental study, we will only consider urban parks as a representative NBS example for applying this technique.

Finally, we use the Gini coefficient \citep{dixon1987bootstrapping} to globally evaluate the effect that an NBS installation has on the resident population. This coefficient is a real number between 0 and 1, widely recognized as a measure of income inequality, where a lower value indicates a fairer distribution of NBS, while higher values suggest a more inequitable allocation \citep{pu2021fairness}.
As an ex-post metric, the Gini coefficient is computed after solving the optimization problem, using the single cell $(i, j) \in G$ with the corresponding value $f_{ij}$ as a statistical unit.

The tests were performed on an Arch Linux OS with an Intel Core i9-11950H 2.60 GHz CPU and 32 GB of RAM. The source code was written in Java 17, with some routines implemented in Python 3.12. The MILP was solved using Gurobi 11.0.2. The time limit was set to 1800 seconds.

\section{Experimental results}\label{results}
This section presents the computational results across six Italian cities. {To evaluate the MILP model under varying urban scales and conditions, we created over $7000$ different instances, considering four different sizes, labeled as XS (Extra Small), S (Small), M (Medium), and L (Large).} Details are given in Table \ref{tab:table_info}, where the grid size $W \times H$ and the number of created instances for each city are reported.

\begin{table}[!ht]
\centering
\footnotesize
\begin{tabular}{|c|c|c|c|c|c|c|c|}
  \hline
 \textbf{Size}&
   $W \times H$&
  \textbf{Bologna} &
  \textbf{Catania} &
  \textbf{Milan} &
  \textbf{Naples} &
  \textbf{Rome} &
  \textbf{Turin} \\
  \hline
\textbf{XS}        & $50 \times 50$         & 443 & 150 & 547 & 368 & 2848 & 426 \\
\textbf{S}         & $100 \times 100$       & 135 & 68  & 150 & 106 & 956  & 122 \\
\textbf{M}         & $200 \times 200$      & 44  & 31  & 51  & 37  & 272  & 37  \\
\textbf{L}         & $300 \times 300$       & 23  & 15  & 26  & 19  & 134  & 19  \\
\hline
\multicolumn{2}{|c|}{\textbf{Total}} & 645 & 264 & 774 & 530 & 4210 & 604 \\ \hline
\end{tabular}%

\caption{Summary of the instances reporting, for each size considered, the total number of instances for each city. The sum of all instances is 7027.}
\label{tab:table_info}
\end{table}

Table \ref{tab:table_forbid_pre} shows the average percentage of both forbidden (`F (\%)') and pre-existing cells (`P (\%)') for each city and instance size. The last row reports the average of these values across all instance sizes (`Size'). Note that these values account for cells being forbidden or pre-existing for at least one of the considered NBSs.

\begin{table}[!ht]
\resizebox{\textwidth}{!}{%
\begin{tabular}{|c|cc|cc|cc|cc|cc|cc|}
\hline
\multicolumn{1}{|c|}{\multirow{2}{*}{\textbf{Size}}} &
  \multicolumn{2}{c|}{\textbf{Bologna}} &
  \multicolumn{2}{c|}{\textbf{Catania}} &
  \multicolumn{2}{c|}{\textbf{Milan}} &
  \multicolumn{2}{c|}{\textbf{Naples}} &
  \multicolumn{2}{c|}{\textbf{Rome}} &
  \multicolumn{2}{c|}{\textbf{Turin}} \\ 
\multicolumn{1}{|c|}{}             & F (\%)     & P (\%)    & F (\%)     & P (\%)    & F (\%)     & P (\%)    & F (\%)    & P (\%)   & F (\%)     & P (\%)    & F (\%)    & P (\%)    \\ \hline
\textbf{XS}                                & 24.31 & 7.65 & 30.72 & 7.32 & 39.04 & 11.64 & 39.49 & 9.77 & 23.07 & 7.20 & 44.19 & 9.03 \\
\textbf{S}                                 & 25.62 & 6.79 & 26.06 & 4.91 & 39.32 & 10.84 & 42.69 & 9.44 & 19.65 & 5.63 & 43.86 & 7.93 \\
\textbf{M}                                 & 33.19 & 6.79 & 29.85 & 2.64 & 42.84 & 9.29  & 47.61 & 8.76 & 20.43 & 5.42 & 47.48 & 7.37 \\
\textbf{L}                                 & 36.85 & 4.58 & 22.92 & 1.33 & 47.81 & 7.92  & 62.28 & 6.83 & 24.63 & 5.37 & 49.21 & 5.91 \\ \hline
\multicolumn{1}{|c|}{\textbf{City Avg.}} & 29.99 & 6.45 & 27.39 & 4.05 & 42.25 & 9.92  & 48.02 & 8.70 & 21.95 & 5.90 & 46.18 & 7.56 \\ \hline
\end{tabular}
}
\caption{{Average percentage of forbidden cells  (`F (\%)') and pre-existing NBS cells  (`P (\%)') for each city and instance size (Size). The final row shows the average value for each city.}}
\label{tab:table_forbid_pre}
\end{table}

To select the available budget, we first define the maximum possible budget by considering an ideal installation of the most expensive NBSs in all the cells of $G$. Then, for each instance, the budget is set to a random value between the 30\% and the 50\% of that amount. Evaluating the impact of the budget on the placement of NBSs is beyond the scope of our tests; therefore, it was chosen to be sufficiently large. Moreover, we assign equal weights to all terms in the objective function to ensure that each is treated with equal importance. However, these weights can be adjusted based on specific business needs or priorities.


\begin{table}[!ht]
\centering
\resizebox{\textwidth}{!}{%
\begin{tabular}{|c|cc|cc|cc|cc|cc|cc|}
\hline
\multirow{2}{*}{\textbf{Size}} & \multicolumn{2}{c|}{\textbf{Bologna}} & \multicolumn{2}{c|}{\textbf{Catania}} & \multicolumn{2}{c|}{\textbf{Milan}} & \multicolumn{2}{c|}{\textbf{Naples}} & \multicolumn{2}{c|}{\textbf{Rome}} & \multicolumn{2}{c|}{\textbf{Turin}} \\
              & \textbf{Time} & \textbf{\%Opt} & \textbf{Time} & \textbf{\%Opt} & \textbf{Time} & \textbf{\%Opt} & \textbf{Time} & \textbf{\%Opt} & \textbf{Time} & \textbf{\%Opt} & \textbf{Time} & \textbf{\%Opt} \\
\hline
\textbf{XS}            & 85.90            & 97.29            & 324.27           & 85.33           & 48.29            & 98.54            & 31.51           & 99.46           & 149.47          & 94.80           & 121.13           & 95.07           \\
\textbf{S}             & 363.32           & 87.41            & 750.39           & 66.18           & 296.65           & 89.33            & 124.11          & 96.23           & 762.33          & 66.84           & 306.75           & 89.34           \\
\textbf{M}             & 702.37           & 79.55            & 1099.21          & 51.61           & 546.93           & 76.47            & 231.98          & 97.30           & 1310.57         & 36.40           & 490.11           & 83.78           \\
\textbf{L}             & 989.67           & 60.87            & 1402.47          & 33.33           & 545.95           & 84.62            & 517.34          & 94.74           & 1255.13         & 45.52           & 762.60           & 73.68           \\
\hline
\end{tabular}
}
\caption{Computational performance of the formulation on test instances. The table reports computation time (in seconds) and the percentage of instances solved to optimality for each group size.}
\label{tab:performance_result}
\end{table}

Table \ref{tab:performance_result} shows the computational performance of the proposed formulation across different city and instance sizes. Two metrics are presented for each city and size: ``Time'', representing the average time (in seconds) required to solve each instance group, and ``\%Solved'', indicating the percentage of instances successfully solved to optimality within the time limit for each group size.
As city size increases from XS to L, the time required to solve each instance rises, while the percentage of solved instances generally decreases. This trend highlights the computational challenges associated with larger urban areas, where achieving optimality within the time limit becomes more difficult. 

Looking at Tables \ref{tab:table_forbid_pre} and \ref{tab:performance_result}, we observe a link between the percentage of forbidden cells and the time required to solve the instances. Specifically, for medium and large instances, Catania and Rome have the lowest percentage of forbidden cells (i.e., 22.92\% and 24.63\% for large cases). This results in a larger number of available cells for locating NBSs, thereby increasing the time required to find an optimal solution. Conversely, Naples, with the highest percentage of forbidden cells (62.28\%), exhibits a lower computation time and a higher percentage of instances optimally solved (517.34 seconds and 94.74\% for large instances). Overall, Naples requires the shortest computation time, followed by Milan, and these two cities also have the highest percentage of instances solved to optimality. In contrast, Catania and Rome have the longest computation times and the lowest percentage of optimally solved instances.

\begin{figure}[!ht]
    \centering
    \includegraphics[scale=0.43]{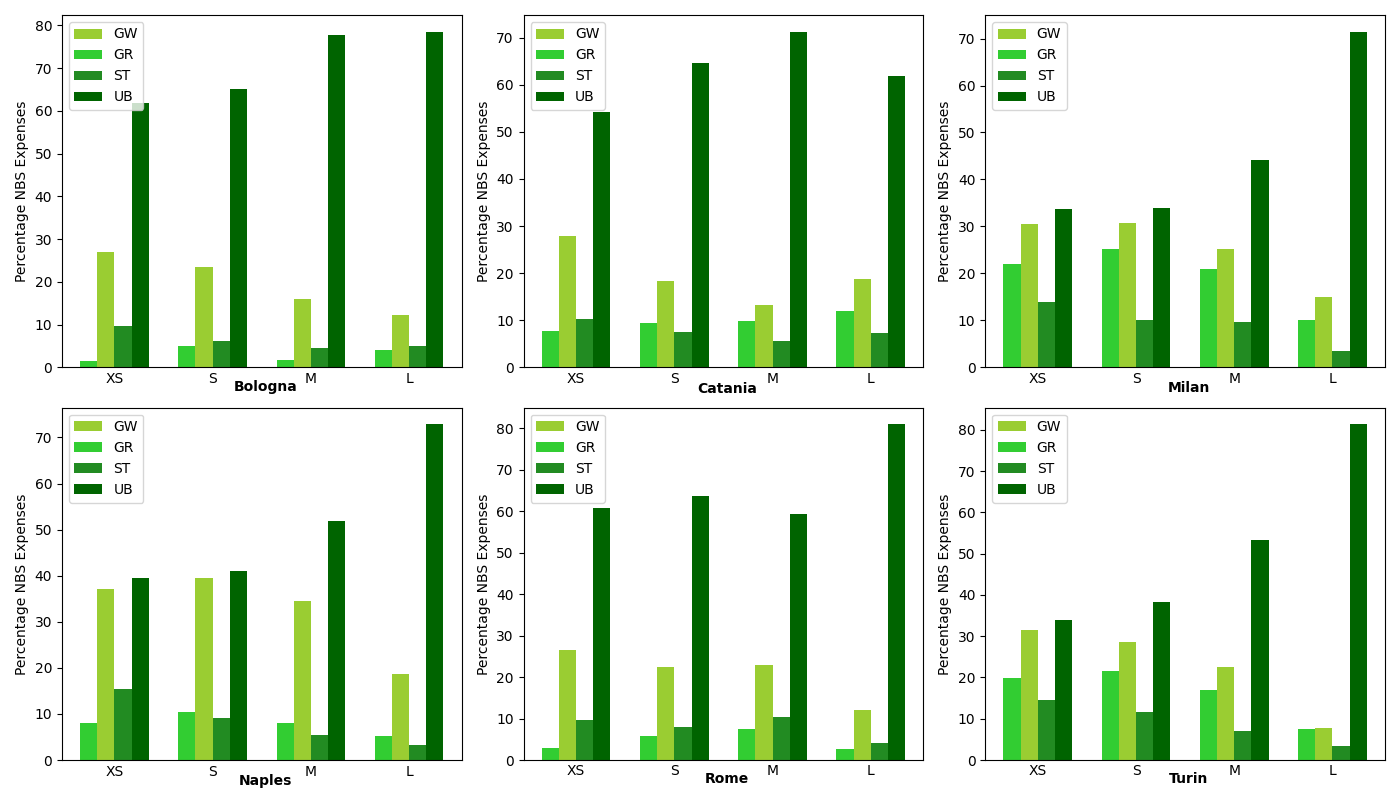}
    \caption{Comparison between the fraction of the budget (in percentage) allocated to each NBS. The reported value represents the average over all the instances of the same size.}
    \label{fig:NBS_percentage_budget_used}
\end{figure}

Figure \ref{fig:NBS_percentage_budget_used} highlights the fraction of the budget (in percentage) allocated to each NBS across the cities, for sizes ranging from XS to L. Urban Parks (UP) dominate budget usage in all cities and sizes, reflecting their greater efficacy in addressing UCs, as supported by the related literature and thus by the larger sizes and higher values of their kernel matrix. Green Walls (GW) and Green Roofs (GR), display a more balanced distribution, although GWs consume a larger portion of the budget compared to GRs, likely due to their greater benefits despite higher installation and maintenance costs. Street Trees (ST) generally require a smaller percentage of the budget, which may reflect a more moderate assessment of their impact or a lower initial cost compared to other NBSs. We note that for ST and GR, the percentage of budget spent decreases as the instance size increases, while for UP the opposite is observed. Clearly, in larger instances, it is possible to locate more clustered NBSs because larger areas are considered.

\begin{figure}[H]
     \centering
     \begin{subfigure}[b]{1\textwidth}
         \centering
         \includegraphics[scale=0.33]{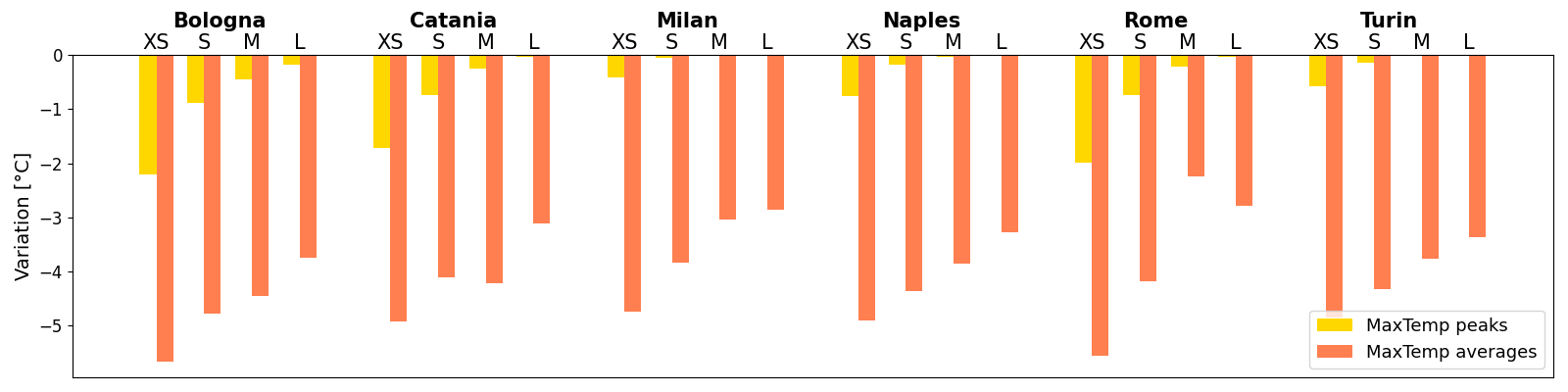}
         \caption{MaxTemp}
         \label{fig:MAX Temp}
    \end{subfigure}
    \\
    \begin{subfigure}[b]{1\textwidth}
         \centering
         \includegraphics[scale=0.33]{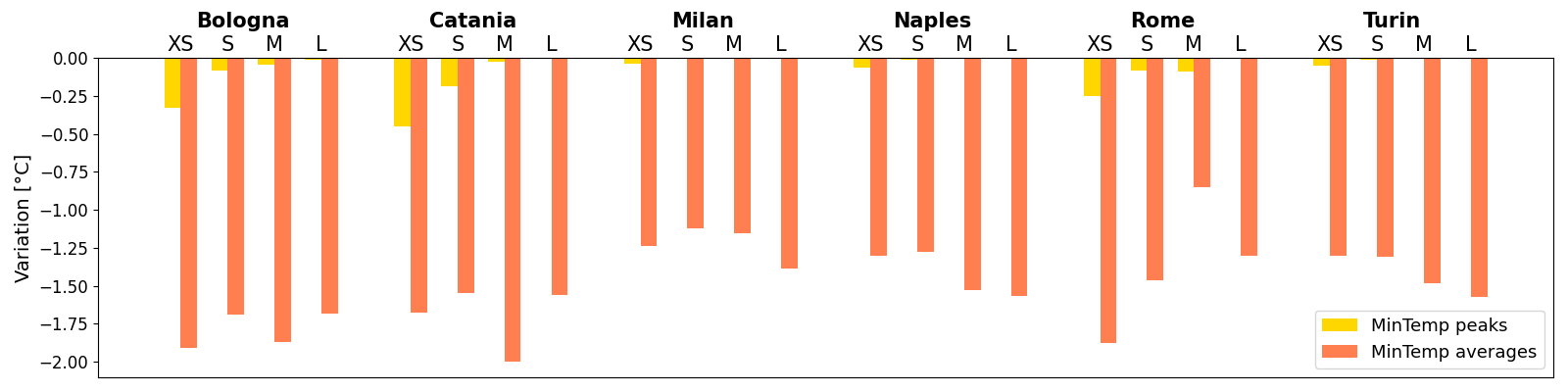}
         \caption{MinTemp}
         \label{fig:MIN Temp}
    \end{subfigure}
    \\
    \begin{subfigure}[b]{1\textwidth}
         \centering
         \includegraphics[scale=0.33]{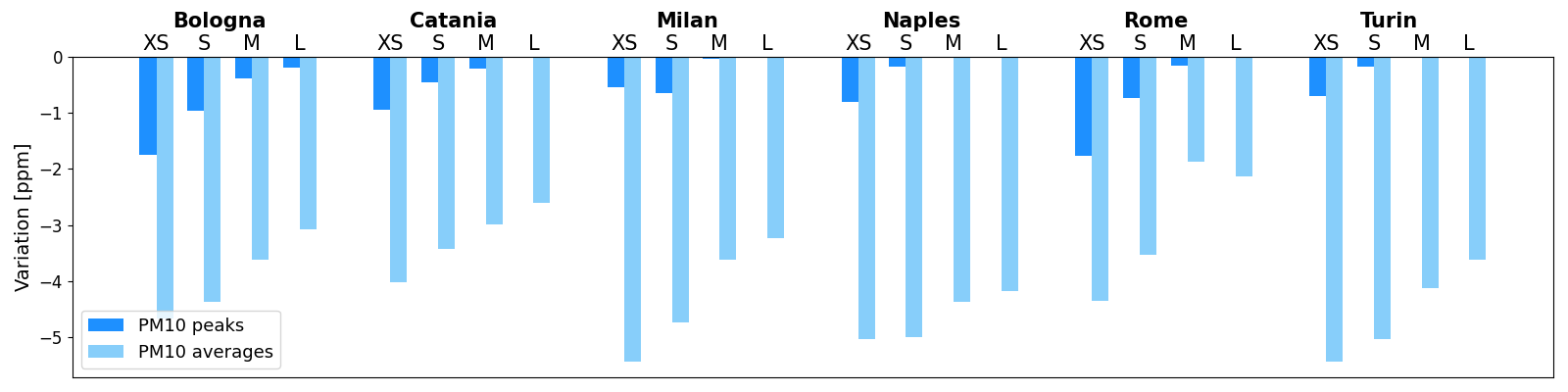}
         \caption{$\textrm{PM}_{10}$}
         \label{fig:PM10}
     \end{subfigure}
     \\
    \begin{subfigure}[b]{1\textwidth}
         \centering
         \includegraphics[scale=0.33]{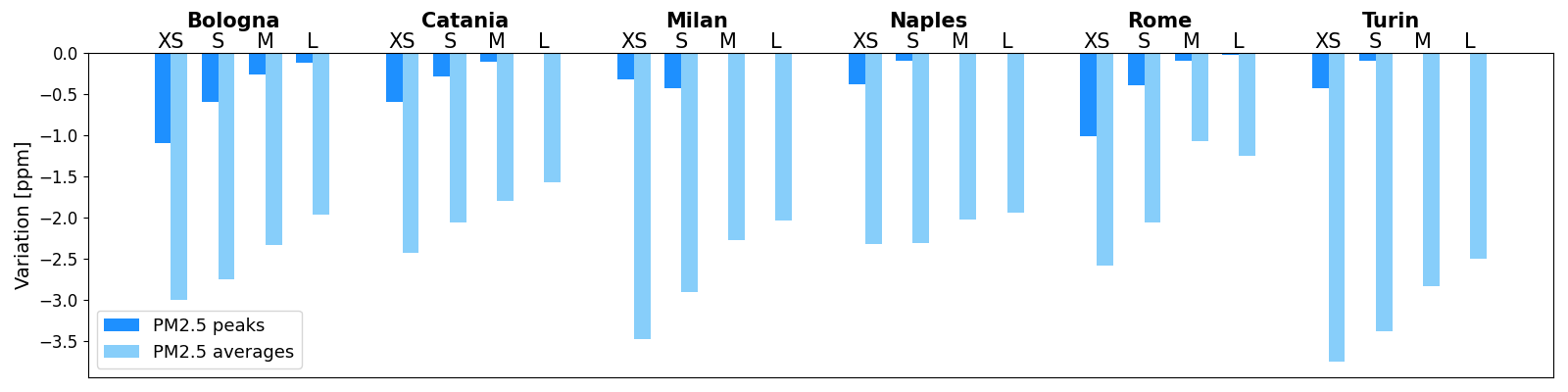}
         \caption{$\textrm{PM}_{2.5}$}
         \label{fig:PM2.5}
     \end{subfigure}
    \caption{Average variations of peaks and average values of four Urban Challenges, respectively (a) MaxTemp variation, (b) MinTemp variation, (c) $\textrm{PM}_{10}$ variation, and (d) $\textrm{PM}_{2.5}$ variation. The variations correspond to reductions of the relative unit (i.e., [°C] for temperature and [ppm] for $\textrm{PM}_{10}$, and $\textrm{PM}_{2.5}$).}
    \label{fig:Average_reduction_of_MAX_and_AVG_UC}
\end{figure}

To highlight how the localized NBSs improve air quality and reduce heat islands, Figure \ref{fig:Average_reduction_of_MAX_and_AVG_UC} presents four sub-figures, each for a specific UC measure. In particular, the sub-figures contain respectively: (a) the maximum temperature (MaxTemp) variation; (b) minimum temperature (MinTemp) variation; (c) $\textrm{PM}_{10}$ variation; (d) $\textrm{PM}_{2.5}$ variation. Each sub-figure, for each city and instance size, shows the average variation over the whole instance and the reduction of the peaks (namely, the maximum variation value of the UC measure over the entire instance). These values are the average across all the instances of the same size.

Note that the model's objective function considers the reduction of both the peaks and the average value of the UC measure of interest. While the model achieves a more substantial reduction in average values across all UCs, the reduction in peak values is smaller. This is because some peaks occur in areas where NBSs cannot be implemented, such as densely built-up clusters with no accessible space. Despite this limitation, the strategic placement of NBSs demonstrates clear benefits, as evidenced by the significant reduction in the average UC measures. We can generally affirm that the average reduction decreases as the instance size increases. However, there are a few cases where the opposite happens (i.e., the minimum temperature for Milan, Naples, and Turin). 

\begin{figure}[H]
    \centering
    \includegraphics[width=\textwidth]{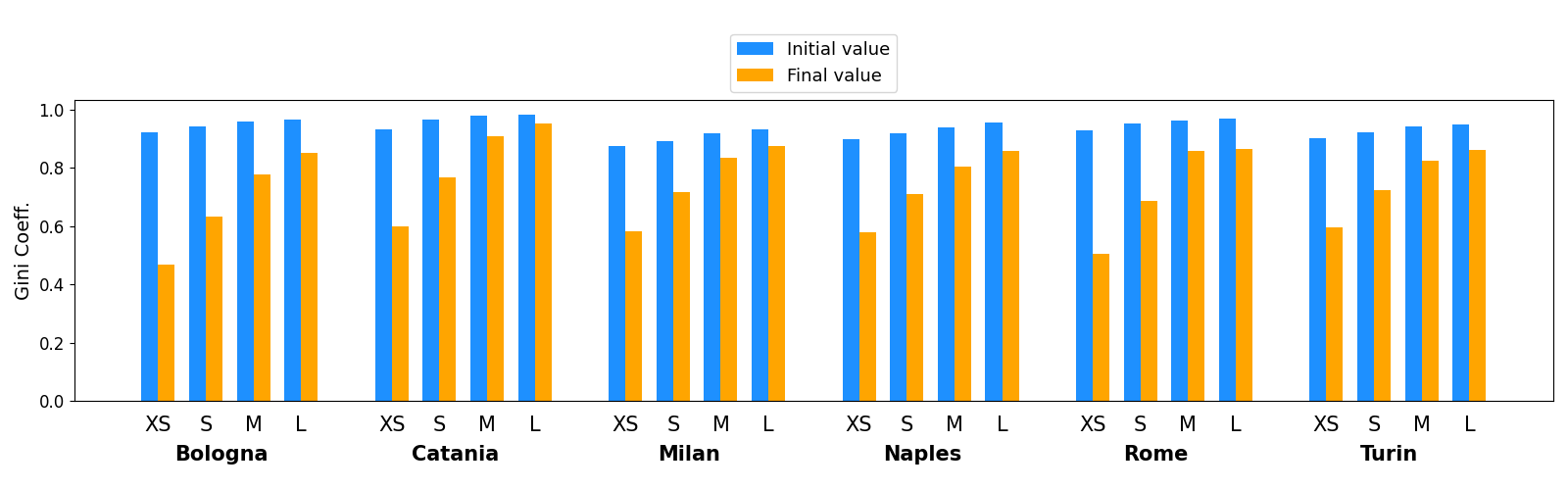}
    \caption{Assessment of the fairness in the placement of NBSs using the Gini coefficient. The figure displays the initial and the final Gini coefficient values for different cities by varying instance sizes.}
    \label{fig:GINI_final_VS_init}
\end{figure}

Fairness in NBS accessibility is evaluated using the Gini coefficient. Figure \ref{fig:GINI_final_VS_init} presents the average Gini coefficient across all instances for each city and size. The Gini coefficient is computed for two scenarios: the initial state (prior to the localization of NBS, considering only existing green spaces) and the final state (after the localization of NBS based on the model's best solution). Across all instances, a decrease in the Gini coefficient is observed, representing enhanced equity in accessing NBS's benefits. Smaller instances exhibit a significant improvement in accessibility. In contrast, larger instances, which cover broader areas, highlight the ability of the proposed approach to account for accessibility by localizing NBSs while considering the distribution of the population.

\begin{figure}[H]
     \centering
     \begin{subfigure}[b]{0.45\textwidth}
         \centering
         \includegraphics[width=\textwidth]{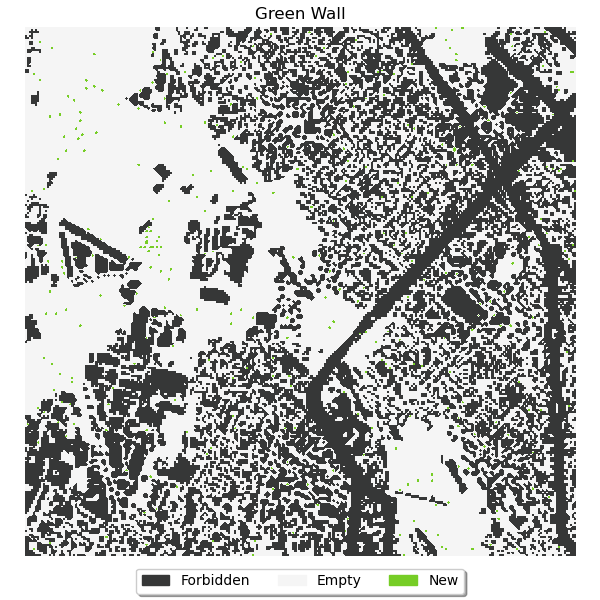}
         \caption{}
         \label{fig:Milano_GW}
     \end{subfigure}
    \begin{subfigure}[b]{0.45\textwidth}
         \centering
         \includegraphics[width=\textwidth]{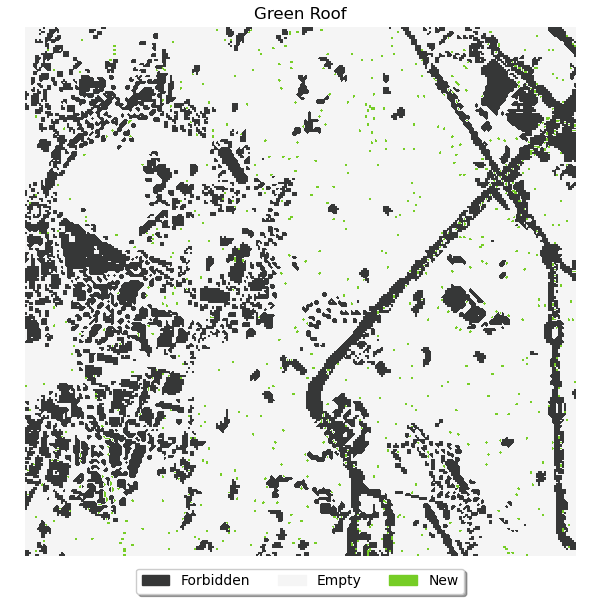}
         \caption{}
         \label{fig:Milano_GR}
     \end{subfigure}
          \\
    \begin{subfigure}[b]{0.45\textwidth}
         \centering
         \includegraphics[width=\textwidth]{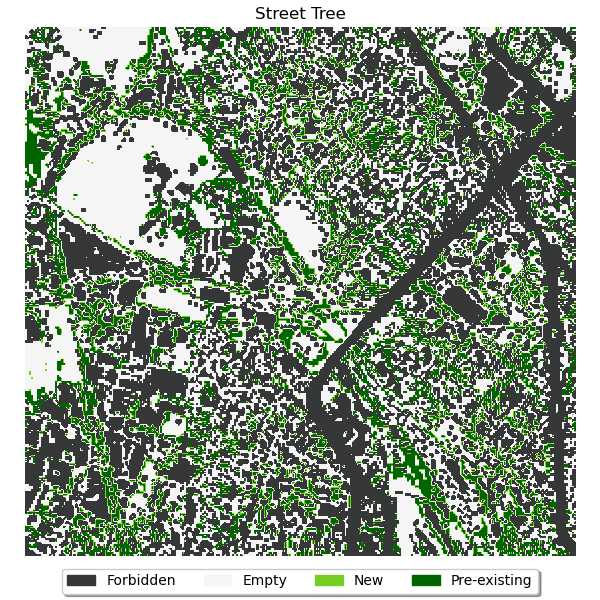}
         \caption{}
         \label{fig:Milano_ST}
     \end{subfigure}
    \begin{subfigure}[b]{0.45\textwidth}
         \centering
         \includegraphics[width=\textwidth]{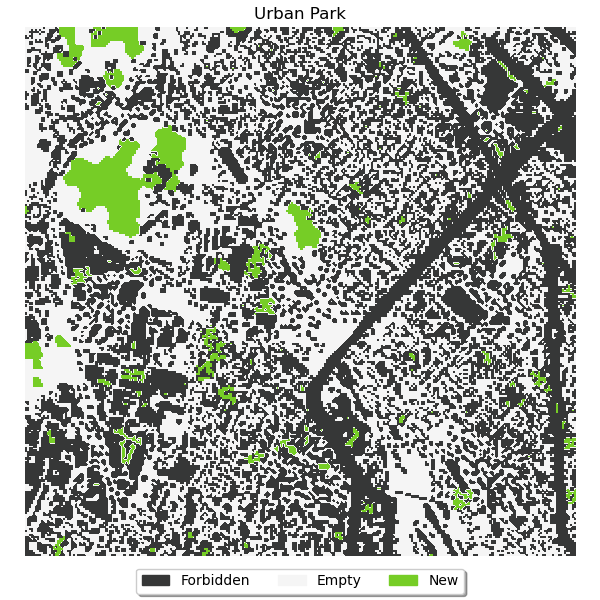}
         \caption{}
         \label{fig:Milano_UP}
     \end{subfigure}
        \caption{NBSs localized in instance L\_10 of MIlan as a result of the MILP model. Each subfigure depicts the distribution of pre-existing and newly installed NBS types (respectively in dark green and light green), as well as areas where localization is not feasible (in dark grey). White areas represent candidate sites not selected for NBS installation by the MILP model.}
        \label{fig:Milano_NBS_Localized}
\end{figure}

\begin{figure}[H]
     \centering
     \begin{subfigure}[b]{0.45\textwidth}
         \centering
         \includegraphics[scale=0.45]{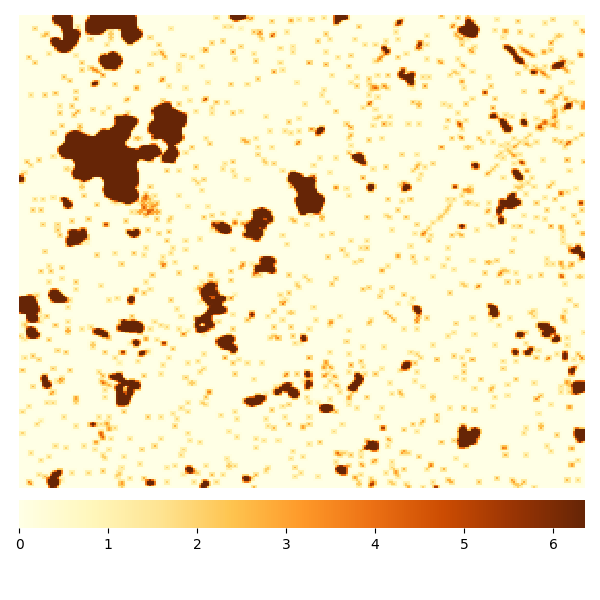}
         \caption{TempMax}
         \label{fig:Milano_TempMax}
     \end{subfigure}
    \begin{subfigure}[b]{0.45\textwidth}
         \centering
         \includegraphics[scale=0.45]{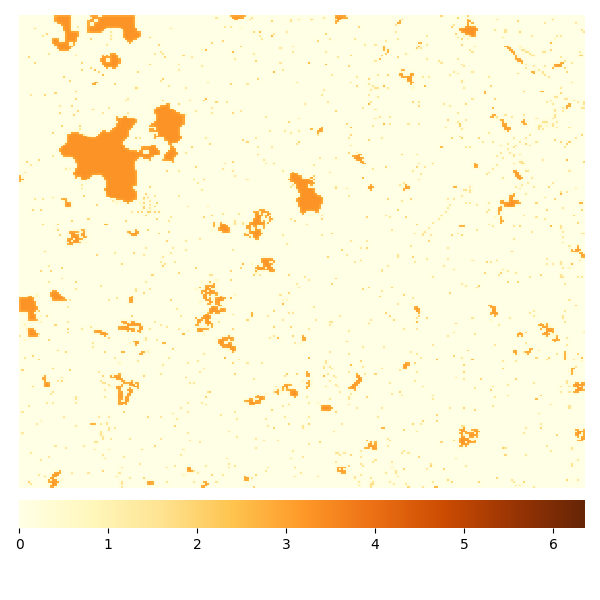}
         \caption{TempMin}
         \label{fig:Milano_TempMin}
     \end{subfigure}
          \\
    \begin{subfigure}[b]{0.45\textwidth}
         \centering
         \includegraphics[scale=0.45]{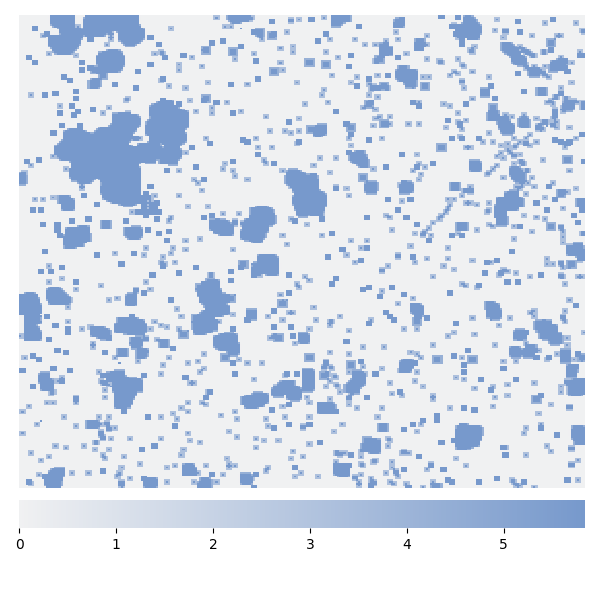}
         \caption{$\textrm{PM}_{10}$}
         \label{fig:Milano_Pm10}
     \end{subfigure}
    \begin{subfigure}[b]{0.45\textwidth}
         \centering
         \includegraphics[scale=0.45]{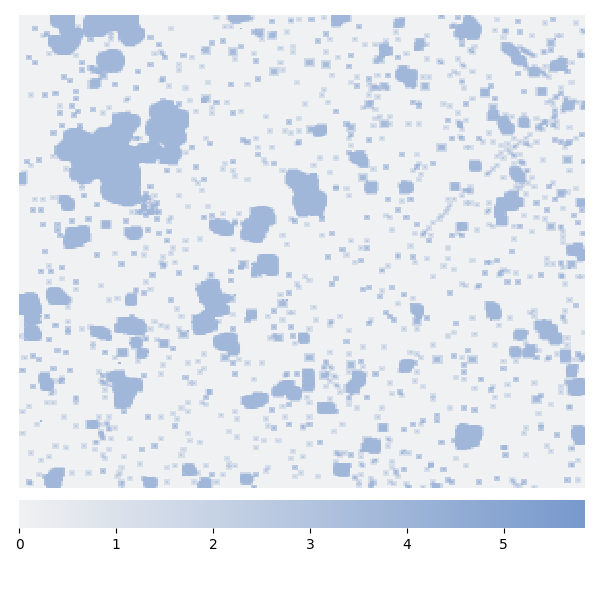}
         \caption{$\textrm{PM}_{2.5}$}
         \label{fig:Milano_Pm2}
     \end{subfigure}
        \caption{Heatmap for each UC measure (TempMax, TempMin, $\textrm{PM}_{10}$, $\textrm{PM}_{2.5}$) showing the difference between the initial measures and the final values obtained after NBSs localization in instance Milan\_L\_10. Darker colors correspond to larger variations.}
        \label{fig:HeatMapMilano}
\end{figure}

To clearly illustrate an example of the output of the optimal NBS planning phase, Figure \ref{fig:Milano_NBS_Localized} shows the landscape of an instance of size L for the city of Milan. Each subfigure presents the distribution of both existing and newly installed NBS types after solving the MILP model (respectively in dark green and light green). Additionally, the figure indicates areas where NBS placement is not feasible and potential sites not selected by the model for NBS installations (in dark grey). To complement this information, Figure \ref{fig:HeatMapMilano} provides a heatmap for each addressed UC measure on the same instance. These heatmaps clearly demonstrate the effect (reduction of all considered quantities) of the installed NBS on the urban landscape. Note that the UC measure reported in each heatmap considers the contribution of all the NBSs installed. A further example for Naples is reported in \ref{app_res}.

\section{Conclusions and perspective outlooks}\label{conclusions}
Recognizing the significance of OR in urban planning, this paper presents a flexible mathematical formulation for locating multiple NBSs while accounting for real-world constraints. We propose a MILP model to optimize the placement of various NBSs in urban environments, addressing key aspects such as ground temperature reduction, air quality improvement, and enhanced accessibility to green spaces for the population. The effectiveness of our methodology is demonstrated through real-world case studies from Italian cities, chosen for their diverse urban characteristics, climatic conditions, and geographic contexts. To the best of our knowledge, this work introduces the first unified OR-based framework for optimizing NBS placement.

\subsection{Managerial Insights}
The proposed framework can provide practical guidance for urban planners and city managers aiming to improve sustainability and resilience through the strategic deployment of NBS. By leveraging the optimization model, city managers can explore and simulate various scenarios by adjusting managerial parameters \citep{TRUONG2016856}, such as budget constraints and objective function weights, which influence the model's outcomes. Thus, city managers can:

\begin{itemize}
    \item \textit{Make Data-Informed Decisions}: Integrating environmental, social, and spatial data, the model empowers decision-makers to prioritize interventions with the highest potential for mitigating urban heat islands and improving air quality.
        
    \item \textit{Recognize High-Impact Opportunities}: Identify and screen out locations that offer limited sustainability benefits and highlight areas with unexploited potential for significant improvements, enabling planners to target locations where NBSs will have the most meaningful environmental and social impacts.
    
    \item \textit{Foster Equitable Interventions}: Balancing environmental benefits across diverse urban areas, promoting equity in NBS implementation, thus ensuring that underserved communities also reap the benefits of urban greening.

    \item \textit{Promote Targeted Incentives:} City managers can use the model’s insights to promote and finance (also through tax discounts) the implementation of green walls and green roofs specifically in areas where these solutions will provide the greatest environmental benefits.
\end{itemize}

\subsection{Research Limitations and Future Directions}
The well-known quote, “\textit{All models are wrong, but some are useful}” by George E. P. Box, was a guiding principle in adjusting the complexity of the method. While no model can perfectly capture reality, its value lies in providing a sufficiently accurate approximation to inform decision-making. To this end, the formulation was intentionally designed as a linear model, leveraging a convolutional-based approach to represent the effects of NBSs. However, it is important to acknowledge that this is a deliberate simplification. Future advancements could involve the integration of non-linear representations of NBS effects or the incorporation of specialized simulation tools.
For instance, integrating advanced simulation tools such as i-Tree ECO \citep{wu2019using} could allow for a more comprehensive assessment of the impacts of NBSs on UCs. These tools, within a simulation-optimization framework, would enable decision-makers to evaluate nonlinear effects and more accurately quantify the benefits of NBSs. Moreover, studies such as \cite{vos2013improving} demonstrate that urban vegetation can produce unintended aerodynamic effects, such as reducing pollutant dispersion, which may counteract its expected benefits. Integrating insights from such research could enhance the model by enabling it to account for trade-offs and localized impacts associated with the implementation of NBSs. 

A further simplification in the model lies in the assumption that NBSs fully occupy the entire grid cell in which they are implemented, without specifying the precise size or configuration of the NBS within that cell. 
Moreover, to enhance the model's accuracy and applicability to real-world scenarios, it can be refined and calibrated using a variety of sources, including historical data, industry reports, future climate projections (such as temperature changes), and expert opinions. These adjustments enable the model to better reflect the complexities of a real-world case study, thereby offering valuable insights into the performance, feasibility, and impact of NBSs under different conditions and contexts.
By addressing these areas, future research can deepen the practical relevance of this optimization framework, offering a robust decision-support tool for sustainable urban planning. Implementing a real-world application would require close collaboration with local administrations to identify feasible areas while accounting for various factors that are beyond the scope of this study. This clarification aims to prevent administrators from raising concerns about possible urban planning or property-related constraints.

\section*{Acknowledgments}
Funding: This work has been supported by (i) CNR DIT.AD106.097 project UISH - Urban Intelligence Science Hub, funded by the PON-METRO Asse 1 Agenda Digitale [CUP: B51B1200430001]; (ii) CNR DIT 
DIT.AD016.094.005 project CTEMT - Casa delle Tecnologie Emergenti di Matera, funded by the PON-METRO Asse 1 Agenda Digitale [CUP: I14E20000020001]; (iii) PNRR MUR project PE0000013-FAIR - Future Artificial Intelligence Research [CUP: B53C22003630006].

\section*{Declaration of interest}
The authors declare that they have no known competing financial interests or personal relationships that could have appeared to influence the work reported in this paper.

\section*{CRediT author statement}
Davide Donato Russo, Diego Maria Pinto, Antonio M. Sudoso: Conceptualization; Methodology; Validation; Formal analysis; Investigation; Data Curation; Writing - Original Draft; Writing - Review \& Editing; Software.

\bibliography{bibliography}

\appendix
\section{}
\subsection{Instance creation}\label{appendix:instance_details}
\noindent
To create the set of benchmark instances, we carried out the following pre-processing steps:

\begin{enumerate}
    \item \textit{Area Definition}: we extracted from each source (i.e., land-use and land-cover data, EEA dataset, temperature dataset, and ISTAT dataset) a sub-matrix containing only the information of the city considered, using the Italian national boundary provided by ISTAT.
    
    \item \textit{Data Reshaping}: all matrices were re-shaped in size to match the target spatial resolution. This becomes necessary because of the different resolutions of the retrieved data. Indeed, when needed, we interpolated the data to fill the missing values in an up-sampling procedure. 
    
    \item \textit{Map Splitting}: to generate multiple instances from the map of a single city, we explored the option of dividing each map into several sub-maps. Thus, given a map $G$, we split it into sub-maps of the following sizes $W \times H$: $50\times50$, $100\times100$, $200\times200$, $300\times300$, and labeled XS, S, M, L, respectively.
\end{enumerate}

Finally, for each generated sub-map, we collected the information to generate an instance and export it as a JSON file. We refer to each instance as $\texttt{CITY\_SIZE\_ID}$, where $\texttt{CITY}$ is the city name, $\texttt{SIZE}$ is the instance label, and $\texttt{ID}$ is a progressive identifier. The instances contains the following information:
\begin{itemize}
    \item $W$ and $H$ as the width and height of the grid, respectively;
    \item the list of NBSs to install, and UCs to address;
    \item for each NBS, the list of cells of the grid that are forbidden for the localization of the NBS, given by land-use information;
   \item for each NBS, the list of cells of the grid where that specific NBS is already located (retrieved from land-use dataset);
    \item for each UC, a matrix describing the current observed measure.
\end{itemize}

\subsection{Kernel selection}\label{appendix:kernel_details}
\noindent
NSBs that deal with green infrastructures have demonstrated improvements in air quality and temperature regulation, in addition to improving human health outcomes associated with air pollution and extreme heat \citep{anderson2022lowering, nowak2018air}. For the considered NBS, we report in Table \ref{tab:NBS-Impact-temp} the reduction of the surface temperature (in \textit{°C}), along with relevant literature sources. 
\begin{table}[H]
    \centering
    \footnotesize
    \begin{tabular}{|c|c|c|}
    \hline
        \textbf{NBS} &  \textbf{Surface Temp. (°C)} & \textbf{References} \\
        \hline
         Green Wall   & 2.7 & \cite{oquendo2022systematic,susca2022effect}  \\
         Green Roof   & 2.0 & \cite{jamei2021review}                        \\
         Street Trees & 1.3 & \cite{segura2022street, gillner2015role}      \\
         Urban Park   & 3.5 & \cite{skoulika2014thermal}                    \\
         \hline
    \end{tabular}
    \caption{Surface temperature reduction in \textit{°C} of the considered NBSs.}
    \label{tab:NBS-Impact-temp}
\end{table}

We considered the work of \cite{di2023cost} to retrieve an NBS performance score related to air pollutants, representing the capacity of the NBS to provide an impact on $\textrm{PM}_{10}$ and $\textrm{PM}_{2.5}$ reduction. These values are reported in the second column of Table \ref{tab:NBS-Impact-PM-perc}. We retrieved from  \cite{SRBINOVSKA2021123306} the average percentage of absorption achieved with green walls for $\textrm{PM}_{2.5}$ and $\textrm{PM}_{10}$, which resulted in a maximum reduction of 25\% for $\textrm{PM}_{2.5}$, and of 37\% for $\textrm{PM}_{10}$. To estimate the percentage of absorption for the remaining NBSs, we combined these values with the air score from \cite{di2023cost}. The full set of values is summarized in Table \ref{tab:NBS-Impact-PM-perc} under the name `PM absorption [\%]'.

\begin{table}[!ht]
    \centering
    \footnotesize
    \begin{tabular}{|c|c|c|c|}
    \hline
    \multirow{2}{*}{\textbf{NBS}} & \multirow{2}{*}{\textbf{Air score}} & \multirow{2}{*}{\textbf{$\textrm{PM}_{2.5}$ absorption [\%]}} & \multirow{2}{*}{\textbf{$\textrm{PM}_{10}$ absorption [\%]}} \\
    & & & \\
    & \citep{di2023cost} & \citep{SRBINOVSKA2021123306} & \cite{SRBINOVSKA2021123306} \\
    \hline
    Green Wall   & 1.00   & 25.00 & 37.00 \\
    Green Roof   & 0.50   & 12.50 & 18.50 \\
    Street Trees & 0.80   & 20.00 & 29.60 \\
    Urban Park   & 1.00   & 25.00 & 37.00 \\
    \hline
    \end{tabular}
    \caption{Estimated percentage of $\textrm{PM}_{2.5}$, $\textrm{PM}_{10}$ absorption taken from the literature. Missing values without references were integrated by considering the air score proposed by \cite{di2023cost}.}
    \label{tab:NBS-Impact-PM-perc}
\end{table}

This percentage of absorption can be converted to values connected to UC measures by considering the average value of PM concentration over the selected use-case. Accordingly, the ranges of the concentration (in $\mu g / m^3$) are the following: $\textrm{PM}_{2.5} \in [5.98, 34.24]$, $\textrm{PM}_{10} \in [1.79, 67.93]$. Thus, Table \ref{tab:NBS-Impact-PM-val} presents the results of applying the percentage reported in Table \ref{tab:NBS-Impact-PM-perc} to the average value of the ranges of pollutant's concentration.

\begin{table}[H]
    \centering
    \footnotesize
    \begin{tabular}{|c|c|c|c|}
    \hline
        \textbf{NBS} & \textbf{$\textrm{PM}_{2.5}$ absorption [$\mu g / m^3$]}  &\textbf{ $\textrm{PM}_{10}$ absorption [$\mu g / m^3$]} \\
        \hline
         Green Wall     & 5.03 & 12.90 \\
         Green Roof     & 2.51 & 6.45  \\
         Street Trees   & 4.02 & 10.32 \\
         Urban Park     & 5.03 & 12.90 \\
         \hline
    \end{tabular}
    \caption{Estimated value of $\textrm{PM}_{2.5}$, $\textrm{PM}_{10}$ absorption in $\mu g / m^3$ .}
    \label{tab:NBS-Impact-PM-val}
\end{table}

We gather NBS data to populate the kernel matrices as follows. For the maximum temperature (TempMax), the central value of the kernel matrix is set equal to the value reported in Table \ref{tab:NBS-Impact-temp}. Then, once the kernel size is defined, the remaining entries in the matrix are linearly reduced, starting from the central highest value and gradually decreasing to values ($> 0$) as the distance from the center increases. We could not retrieve a suitable estimation for the minimum temperature (TempMin). However, according to \cite{DIPIRRO20221}, the minimum and maximum temperatures differ by approximately 70\%. Based on this, we set the kernel entries for TempMin by proportionally reducing the values used for TempMax by this percentage. Kernels for pollutants were generated from the values presented in Table \ref{tab:NBS-Impact-PM-val} as well as for the temperature kernels generated from values presented in Table \ref{tab:NBS-Impact-temp}. The remaining set of kernels, one for each NBS, evaluates the fairness in NBS localization. The entries within these kernels are manually selected by considering the NBS potential to have a greater impact on the population, ensuring that areas with higher population density receive adequate attention.

\section{Results} \label{app_res}
\begin{figure}[!ht]
     \centering
     \begin{subfigure}[b]{0.45\textwidth}
         \centering
         \includegraphics[width=\textwidth]{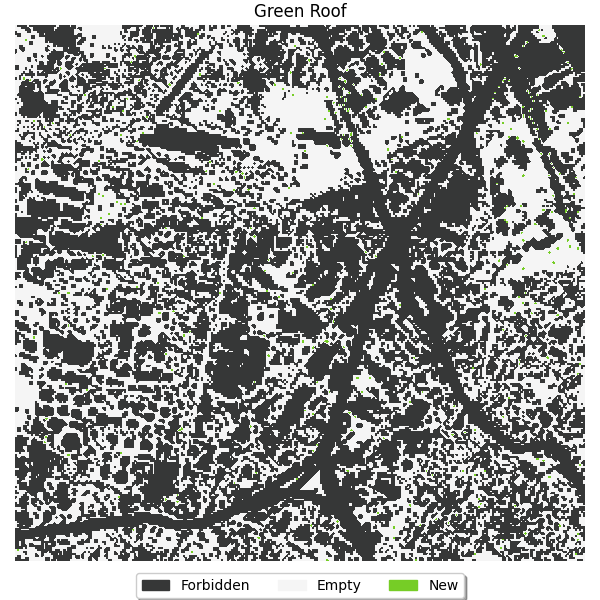}
         \caption{}
         \label{fig:Naples_GW}
     \end{subfigure}
    \begin{subfigure}[b]{0.45\textwidth}
         \centering
         \includegraphics[width=\textwidth]{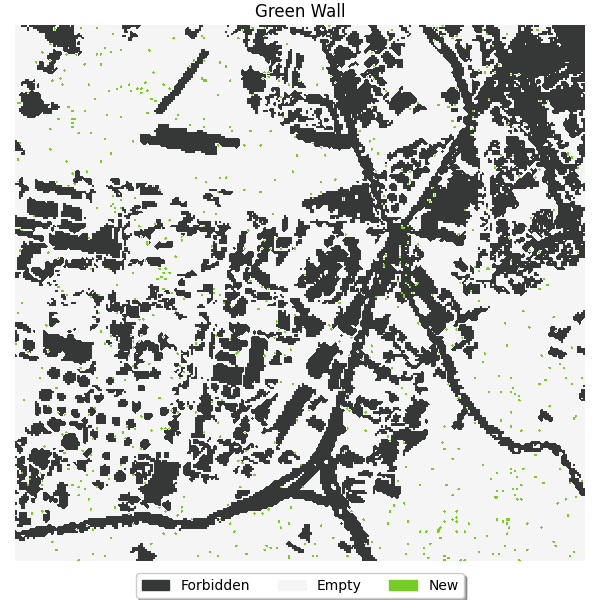}
         \caption{}
         \label{fig:Naples_GR}
     \end{subfigure}
          \\
    \begin{subfigure}[b]{0.45\textwidth}
         \centering
         \includegraphics[width=\textwidth]{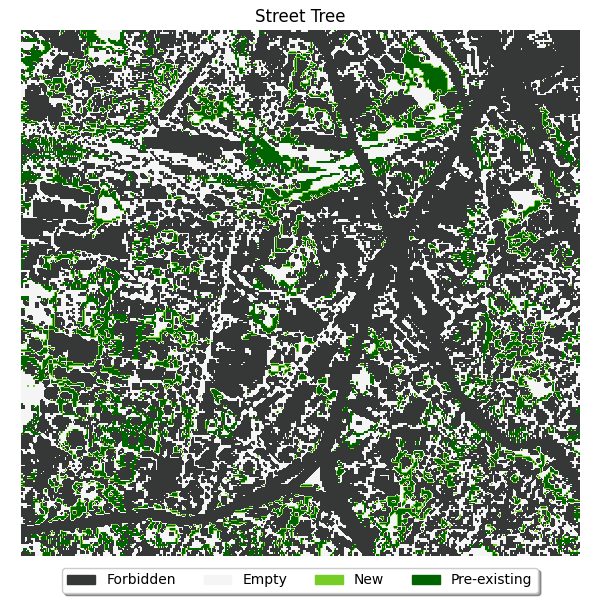}
         \caption{}
         \label{fig:Naples_ST}
     \end{subfigure}
    \begin{subfigure}[b]{0.45\textwidth}
         \centering
         \includegraphics[width=\textwidth]{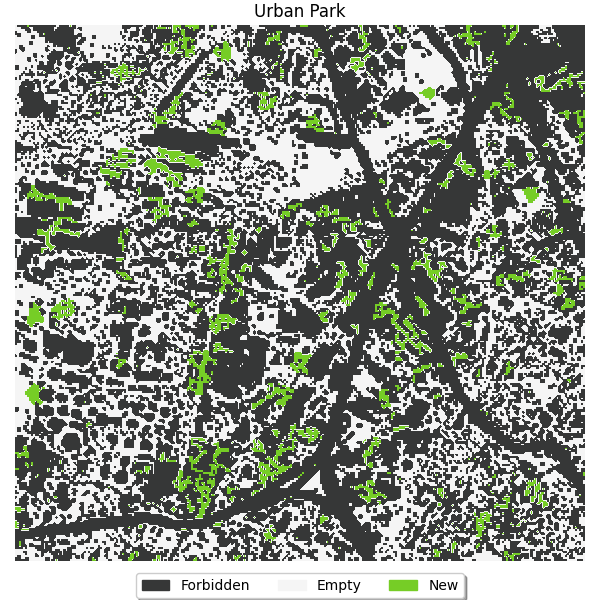}
         \caption{}
         \label{fig:Naples_UP}
     \end{subfigure}
        \caption{NBSs localized in instance L\_18 of Naples as a result of the MILP model. Each subfigure depicts the distribution of pre-existing and newly installed NBS types, as well as areas where localization is not feasible. White areas represent candidate sites not selected for NBS installation by the MILP model.}
        \label{fig:Naples_NBS_Localized}
\end{figure}

\begin{figure}[!ht]
     \centering
     \begin{subfigure}[b]{0.45\textwidth}
         \centering
         \includegraphics[scale=0.45]{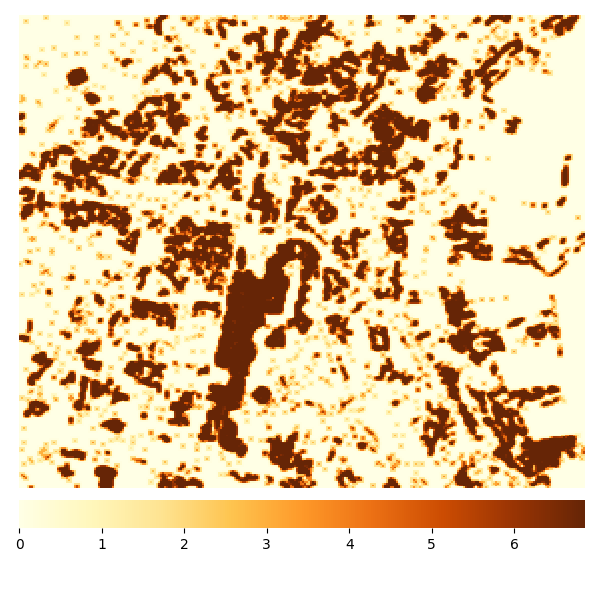}
         \caption{TempMax}
         \label{fig:Naples_TempMax}
     \end{subfigure}
    \begin{subfigure}[b]{0.45\textwidth}
         \centering
         \includegraphics[scale=0.45]{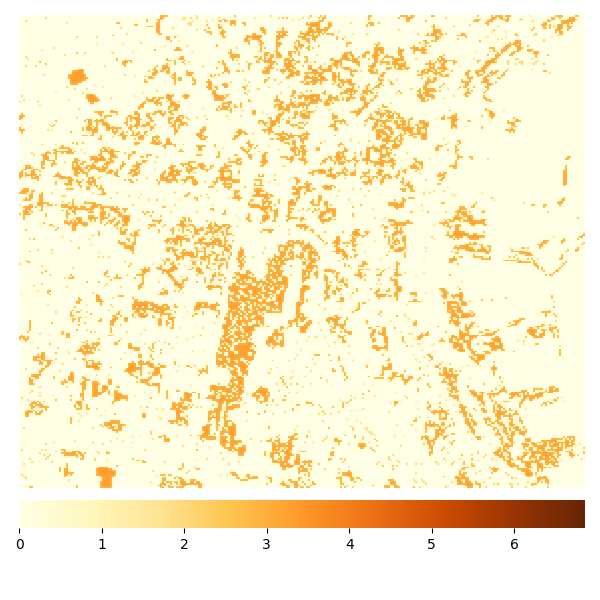}
         \caption{TempMin}
         \label{fig:Naples_TempMin}
     \end{subfigure}
          \\
    \begin{subfigure}[b]{0.45\textwidth}
         \centering
         \includegraphics[scale=0.45]{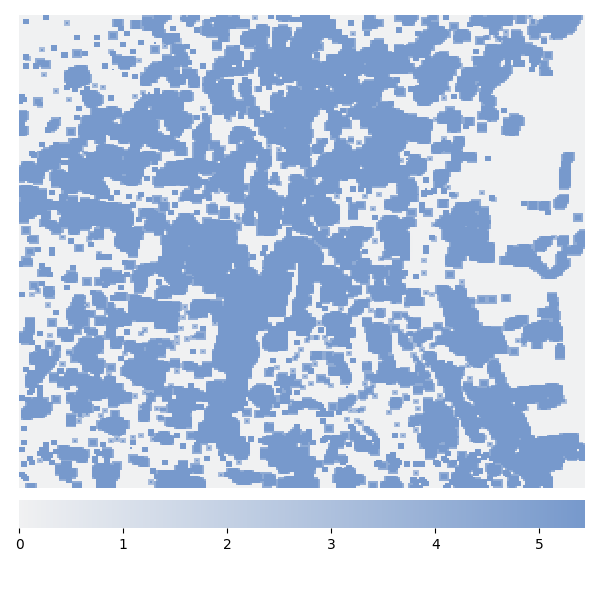}
         \caption{$\textrm{PM}_{10}$}
         \label{fig:Naples_Pm10}
     \end{subfigure}
    \begin{subfigure}[b]{0.45\textwidth}
         \centering
         \includegraphics[scale=0.45]{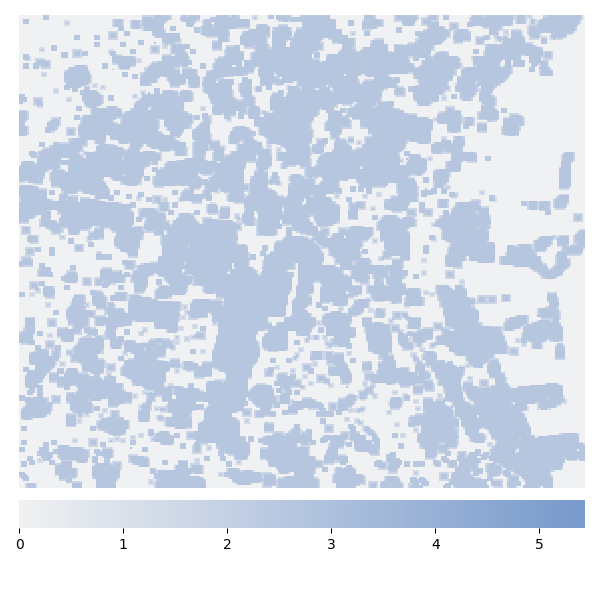}
         \caption{$\textrm{PM}_{2.5}$}
         \label{fig:Naples_Pm2}
     \end{subfigure}
        \caption{Heatmap for each UC measure (TempMax, TempMin, $\textrm{PM}_{10}$, $\textrm{PM}_{2.5}$) showing the difference between the initial measures and the final values obtained after NBSs localization in instance Naples\_L\_18. Darker colors correspond to larger variations.}
        \label{fig:HeatMapNaples}
\end{figure}

A further example of the outputs of the optimal NBS planning phase for one more instances is reported in this Appendix. Figure \ref{fig:Naples_NBS_Localized} shows the landscape of the instance of size L for the city of Naples (i.e. Naples\_L\_18). Each subfigure presents the distribution of both existing and newly installed NBS types after solving the MILP model (respectively in dark green and light green). The figures indicate areas where NBS placement is not feasible and potential sites not selected by the model for NBS installations (in dark grey). To highlight the effect of NBS installations, Figure \ref{fig:HeatMapNaples} provides a heatmap for each addressed UC measure.

\end{document}